\documentclass{amsart}
\usepackage{latexsym,amsmath,amssymb,amsfonts,amscd,graphics}
\usepackage[mathscr]{eucal}
\usepackage[all]{xypic}
\usepackage{hyperref}

\newtheorem{theorem}{Theorem}[section]
\newtheorem{lemma}[theorem]{Lemma}
\newtheorem{proposition}[theorem]{Proposition}
\newtheorem{corollary}[theorem]{Corollary}

\theoremstyle{definition}
\newtheorem{definition}[theorem]{Definition}
\newtheorem{notation}[theorem]{Notation}

\theoremstyle{remark}
\newtheorem{remark}[theorem]{Remark}

\numberwithin{equation}{section}

\newcommand{\bP}{\mathbb{P}}

\newcommand{\bC}{\mathbb{C}}

\newcommand{\calM}{\mathcal{M}}
\newcommand{\calO}{\mathcal{O}}

\newcommand{\calI}{\mathcal{I}}

\newcommand{\Sym}{\mathrm{Sym}}
\newcommand{\PGL}{\mathrm{PGL}}
\newcommand{\SL}{\mathrm{SL}}

\newcommand{\Sing}{\mathrm{Sing}}

\newcommand{\Sec}{\mathrm{Sec}}

\newcommand{\gquot}{/\!\!/}
 
\author{Radu Laza}
\title[Moduli of cubic $4$-folds]{The moduli space of cubic fourfolds} 
\address{University of Michigan \\
1832 East Hall \\
Ann Arbor, MI 48109}
\email{rlaza@umich.edu}
\begin{document}
\bibliographystyle{amsplain}
\begin{abstract}
We describe the GIT compactification of the moduli space of cubic fourfolds (cubic hypersurfaces in the five dimensional projective space), with a special emphasis on the role played by  singularities. Our main result is that a cubic fourfold with only isolated simple (A-D-E) singularities is GIT stable. Conversely, with some minor exceptions, the stability for cubic fourfolds is characterized by this condition. 
\end{abstract}
\maketitle

\section{Introduction}\label{sectintro}
One of the classical results of algebraic geometry is that the moduli space of elliptic curves can be constructed both as a GIT quotient and as the quotient of the upper half plane by the modular group. Several similar examples, where a dual construction for the moduli space exists, were known classically (e.g. the  low degree $K3$ surfaces), and   a few more were discovered recently (e.g. the moduli spaces of cubic  surfaces and cubic threefolds -- see \cite{allcock2}, \cite{allcock3fold} and \cite{looijengaswierstra}). The purpose of this and of a subsequent paper is to discuss a new example in this vein, namely, we analyze the moduli space of cubic fourfolds.

\smallskip
 
The moduli space of cubic fourfolds came to attention recently as a key ingredient in the  Allcock--Carlson--Toledo \cite{allcock3fold} and Looijenga--Swierstra \cite{looijengaswierstra} construction of a uniformization of the moduli  space of cubic threefolds by a complex ball.  One important reason for this is that the period map for cubic fourfolds behaves quite nicely, very similarly to the period map for K3 surfaces. There are a series of results in this direction (e.g. Beauville--Donagi \cite{beauvilledonagi}) culminating with the proof of the global Torelli theorem by Voisin \cite{voisin} and the work of Hassett \cite{hassett0} on rationality.  However, in contrast to the $K3$ case,  surjectivity type results for cubic fourfolds (i.e. the characterization of the image of the period map) are not known. 

\smallskip

Following the example of low degree K3 surfaces (Shah  \cite{shah,shah4}) and cubic threefolds (see \cite{allcock1},  \cite{yokoyama}, \cite{allcock3fold} and \cite{looijengaswierstra}), we attack the  problem mentioned above  by means of geometric invariant theory. Namely, here we discuss the degenerations of cubic fourfolds from a GIT point of view. In subsequent work, these results are combined with a  monodromy analysis to get a good understanding of the image of the period map for cubic fourfolds. The main idea behind this approach is that, quite generally, GIT provides a quick way of compactifying a moduli space. The resulting compactification can then be used as a starting point for more refined constructions.

\smallskip

The main problem, however, is that typically the GIT compactification is quite badly behaved. For example, it is quite possible that  nice objects for the moduli problem under consideration are left out of the GIT compactification. Our first main result is that this type of phenomena does not occur for  cubic fourfolds.  Specifically, the following theorem essentially identifies the  stable cubic fourfolds as those having simple singularities. 

\smallskip

\begin{theorem}\label{mainthm1}
A cubic fourfold $Y$ is not stable if and only if one of the following conditions holds:
\begin{itemize}
\item[i)] $Y$ is singular along a curve $C$ spanning a linear subspace of dimension at most $3$ of $\bP^5$; 
\item[ii)] $Y$ contains a singularity that deforms to a singularity of class $\widetilde{E_r}$ (for $r=6,7,8$). 
\end{itemize}
In particular, if $Y$ is a cubic fourfold with isolated singularities, then $Y$ is stable if and only if $Y$ has at worst simple singularities. 
\end{theorem}

\smallskip

In particular, Theorem \ref{mainthm1}  allows us to speak about the moduli space of cubic fourfolds $\calM$ with at worst simple singularities. This is extremely relevant in the context of analyzing the period map for cubic fourfolds. We recall that the simple singularities in even dimensions are characterized by the fact that they give finite monodromy. Therefore, the space $\calM$ is the natural space where the period map (defined a priori only for smooth cubic fourfolds) would extend. The fact that the simple singularities give stable points plays a key role in the case of K3 surfaces (see Shah \cite{shah}). It is also implicitly used in the case of cubic threefolds by Allcock--Carlson--Toledo \cite{allcock3fold} and Looijenga--Swierstra \cite{looijengaswierstra}.

\smallskip

Our second main result is the following description of the GIT compactification of the moduli space of cubic fourfolds:

\smallskip

\begin{theorem}\label{mainthm2}
The moduli space  $\calM$ of cubic fourfolds having at worst simple  singularities is compactified by the GIT quotient $\overline{\calM}$ by adding six irreducible boundary components, that we label $\alpha, \dots, \phi$. A semistable cubic fourfold $X$ with minimal orbit corresponding to a generic point in a boundary component has the following geometric property:
\begin{itemize}
\item[Case] $\alpha$: $X$ is singular along a line and a quartic elliptic curve;
\item[Case] $\beta$: $X$ has two $\widetilde{E_8}$ singularities (of the same modulus); 
\item[Case] $\gamma$: $X$ is singular along a conic and has an isolated singularity of $\widetilde{E_7}$; 
\item[Case] $\delta$: $X$ has three $\widetilde{E_6}$ singularities (of the same modulus); 
\item[Case] $\epsilon$: $X$ is singular along a rational normal curve of degree $4$; $X$ is stable;
\item[Case] $\phi$: $X$ is singular along a sextic elliptic curve; $X$ is stable;
\end{itemize}
(see table \ref{tableboundary}). Furthermore, the boundary components $\beta$ and $\epsilon$ are $3$-dimensional and they meet along a surface $\sigma$. The surface $\sigma$ meets the $2$-dimensional boundary components $\gamma$ and $\phi$ along a curve $\tau$. Finally, the curve $\tau$ meets the $1$-dimensional components $\alpha$ and $\delta$ in a point $\zeta$ (see figures \ref{gitboundary} and \ref{gitboundary2}). 
\end{theorem}

\smallskip

The description of the GIT compactification $\overline{\calM}$ might seem  complicated, but we note that there is quite a bit of structure. To start, we note that, as discussed in section \ref{sectcomments}, the GIT computation for cubic fourfolds is closely related to that for cubic threefolds (Allcock \cite{allcock1} and Yokoyama \cite{yokoyama})  and that for plane sextics (Shah \cite{shah}). It follows that  one can essentially reconstruct the cubic fourfold case from these two lower dimensional cases. At a deeper level, in all three cases mentioned here, the structure of the moduli space is dictated by the Hodge theoretical properties of the varieties under consideration. We only lightly touch on this in section \S\ref{secstratification}. Nonetheless, it is quite apparent from our computations that the relationship between the GIT construction and Hodge theoretical construction of the moduli space of cubic fourfolds is very similar to that for low degree $K3$ surfaces (see \cite{shah,shah4} and \cite[\S8]{looijengacompact}).

\smallskip

A few words about the organization of the paper. A standard GIT analysis consists of three steps. The first one is a purely combinatorial one, and consists of identifying certain maximal subsets of monomials. We discuss this step for cubic fourfolds in section \ref{sectgit}.   The next step (corresponding to sections \ref{sectnonstable} and \ref{sectminorbits} in our text)  attaches some geometric meaning to the combinatorial results obtained in the previous step. The results in this step typically describe the stability of hypersurfaces in terms of a ``bad flag'' (see Theorem \ref{thmsstable} and the discussion from \cite[\S4.2]{GIT}). Unfortunately, this geometric interpretation is rather coarse, so one needs to refine these results. Typically, by using some classification of singularities, one can interpret the existence of bad flags in terms of singularities. In our situation, we divide the analysis in two cases: isolated (section \ref{sectisolated})  or non-isolated (section \ref{sectnonisolated}) singularities. The main results here are Theorems \ref{thmsimplesing} and \ref{thmsingularities}. We note that the most delicate aspect is the identification of the boundary strata $\epsilon$ and $\phi$ corresponding to stable cubic fourfolds with non-isolated singularities. Finally, in sections \ref{sectconclusion} and  \ref{sectcomments}, we put everything together and conclude the proofs of the main results and make some further comments on the structure of the GIT compactification $\overline{\calM}$.

\smallskip

After the completion of this work, we've learned that both Allcock \cite{allcockp} and  Yokoyama \cite{yokoyama2} (independently)  have done  a partial analysis of the stability for cubic fourfolds. Our results, however, are more detailed and amenable to the study of the period map for cubic fourfolds.

\smallskip

\subsection{Acknowledgments} While writing this paper, I've benefited from discussions with Igor Dolgachev, Bob Friedman, and Rob Lazarsfeld.   I am grateful for their helpful comments   and support.  I would also like to thank Eduard Looijenga for informing me about his computation of the Baily--Borel compactification and the work of Yokoyama. 

\smallskip

\subsection{Notations and Conventions}\label{sectnotations} We fix the following notations:
\begin{itemize}
\item[-] $\calM_0$: the  moduli space of smooth cubic fourfolds; 
\item[-] $\calM$ : the moduli space of cubic fourfolds with  simple singularities;
\item[-] $\calM^s$ : the moduli space of  stable cubic fourfolds;
\item[-] $\overline{\calM}$ : the GIT compactification of $\calM_0$.
\end{itemize}
Let $W$ ($=H^0(\bP^5,\calO_{\bP^5}(1))$) be a six dimensional vector space. By definition, 
$$\overline{\calM}:=\bP(\Sym^3 W)^{ss}\gquot \SL(W),$$
where the quotient is taken in the sense of GIT \cite{GIT}. $\calM_0$, $\calM$, and $\calM^s$ are open (Zariski) subsets in $\overline{\calM}$. By Theorem \ref{mainthm1}, we have 
 $\calM_0\subset \calM\subset \calM^s\subset \overline{\calM}$. Depending on the context, boundary means either $\overline{\calM}\setminus \calM^s$  or  $\overline{\calM}\setminus \calM$. 
 
 \smallskip
 
 The GIT terminology is that of Mumford \cite{GIT}. For us, {\it unstable} means not semistable, {\it non-stable} is a shorthand for not properly stable, and {\it strictly semistable} means semistable, but not properly stable. As is customary, the one parameter subgroups ($1$-PS) $\lambda$ of $\SL(6)$ that are used in the application of the numerical criterion are  assumed diagonal $t\in \bC^*\xrightarrow{\lambda} \mathrm{diag}(t^{a_0},\dots, t^{a_5})\in \SL(6)$ with the weights satisfying $a_0\ge \dots\ge a_5$ and $a_0+\dots+a_5=0$.   We call such a $1$-PS {\it normalized} and denote it by its weights $(a_0,\dots,a_5)$. Given a monomial $x_0^{i_0}\dots x_5^{i_5}$, its weight with respect to a normalized $1$-PS $\lambda$ of weights $(a_0,\dots,a_5)$ is $i_0\cdot a_0+\dots i_5 \cdot a_5$. We denote by $M_{\le 0}(\lambda)$ (and $M_{<0}(\lambda)$) the set of monomials of degree $3$ which have non-positive (resp. negative) weight with respect to $\lambda$. Similarly, $M_0(\lambda)$ is set of monomials of weight $0$ (i.e. invariant with respect to $\lambda$).

\smallskip

We make the convention that all the boundary strata ($\alpha,\beta,\dots$) that occur in our classification are closed and irreducible. In particular, it makes sense to discuss about a generic point in a boundary stratum. By  an adjacency of strata (represented by an arrow in figures \ref{gitboundary} and \ref{gitboundary2}), we understand an inclusion of strata. The various symbols used in figure \ref{gitboundary} and \ref{gitboundary2} represent the dimensions: $\bullet$, $\circ$, $\Box$, and $\diamond$ mean $1$-, $2$-, $3$-, and $4$-dimensional boundary components respectively. 

\smallskip

We are concerned with the following analytic types of isolated hypersurface singularities: $A_n$ ($n\ge 1$), $D_n$ ($n\ge 4$), $E_r$ ($r=6,7,8$), and $\widetilde{E_r}$ ($r=6,7,8$). The singularities of type $A_n$, $D_n$, and $E_r$ are called {\it simple singularities}. The terminology and notations are those of Arnold \cite{AGV1}, with the  exception that we use  $\widetilde{E_6}$, $\widetilde{E_7}$, $\widetilde{E_8}$ instead of $P_8$, $X_9$, and $J_{10}$ respectively. We also consider two types of non-isolated hypersurface singularities: $A_{\infty}$ (double line) and $D_{\infty}$ (pinch point) given locally in $(\bC^5,0)$  by the equations $(x_2^2+\dots+x_5^2=0)$ and $(x_1x_2^2+\dots+x_5^2=0)$ respectively. For the reader familiar with these notations for surface singularities, we mention that we refer to the fourfold singularities obtained  by suspending the homonymous types of surface singularities.  

\smallskip

By {\it the corank of the hypersurface singularity} given by  $f(x_1,\dots,x_n)$, we understand the number of variables $n$ minus the rank of the Hessian of $f$ (see \cite[Ch. 11]{AGV1}). The corank is a stable  invariant of the singularity. We note  that a hypersurface singularity of corank $1$ is of type $A_n$ (for some $n\ge2$) if isolated, or $A_{\infty}$ otherwise  (see \cite[Ch. 11, 16]{AGV1}). Similarly, a singularity of corank $2$ and third jet $x_1^2x_2$ is of type $D_n$ (for some $n\ge 5$) or $D_{\infty}$.  

\section{Preliminary Study of Stability}\label{sectgit}
The main tool of analyzing the stability for cubic fourfolds is the Hilbert-Mumford numerical criterion \cite[Thm. 2.1]{GIT}:
{\it A cubic form $f$ is stable (resp. semistable) iff $\mu(f,\lambda)>0$ (resp. $\ge 0$) for all $\lambda$ one parameter subgroups of $\SL(6)$}, where $\mu(f,\lambda)$ is the numerical function of Mumford. As is customary, we fix coordinates on $\bP^5$ and assume that all $1$-PS used are normalized as in \S\ref{sectnotations} (equivalently, fix a maximal torus $T$ in $\SL(6)$, and consider only $1$-PS of $T$). In order to understand the stable cubic fourfolds, it suffices to find the cubic  forms $f$ for which there exists a normalized $1$-PS $\lambda$ such that $\mu(f,\lambda)\le 0$ (any other non-stable cubic is projectively equivalent to such an $f$).  By definition (see \cite[pg. 81]{GIT}), $\mu(f,\lambda)$ is the highest weight with respect to $\lambda$ of a monomial occurring with non-zero coefficient in $f$. Thus, the cubic defined by $f$ is unstable if all the monomials in $f$ belong to $M_{\le 0}(\lambda)$ for some $\lambda$. It is then clear that, for the analysis of stability, it is enough to identify the maximal possible subsets $M_{\le 0}(\lambda)$ and to interpret geometrically what it means for $f$ to be a linear combination of monomials from  $M_{\le 0}(\lambda)$.  In our situation, we identify the maximal $M_{\le 0}(\lambda)$ by applying a simple computer program (see remark \ref{remalgo}). The results are given in table \ref{maxsstable} below. We note that  similar computations occur in literature (see esp. \cite{allcock1} and \cite{yokoyama}).

\begin{table}[htb]
\begin{center}
\renewcommand{\arraystretch}{1.25}
\begin{tabular}[2cm]{|c|l|l|c|}
\hline
Case& Weights of $\lambda$& Maximal monomials &Invariant\\
\hline\hline
S1&$( 2, 2,-1,-1,-1,-1)$&  $x_0 x_2^2 $&$\alpha$\\
\hline
S2&$( 1, 1, 1, 1,-2,-2)$& $x_0^2 x_4 $  &$\alpha$\\
\hline
S3&$( 2, 1, 0, 0,-1,-2)$& $x_0 x_2 x_5 $, $x_0 x_4^2 $, $x_1^2 x_5 $, $x_1 x_2 x_4 $, $x_2^3 $ &$\beta$\\
\hline
S4&$( 4, 1, 1,-2,-2,-2)$&   $x_0 x_3^2 $,  $x_1^2 x_3 $&$\gamma$\\
\hline
S5&$( 2, 2, 2,-1,-1,-4)$&  $x_0^2 x_5 $,    $x_0 x_3^2 $&$\gamma$\\
\hline
S6&$( 2, 0, 0, 0,-1,-1)$&  $x_0 x_4^2 $,   $x_1^3 $ &$\delta$\\
\hline
S7&$( 1, 0, 0, 0, 0,-1)$&  $x_0 x_1 x_5 $, $x_1^3 $ \\
\cline{1-3}
S8&$( 1, 1, 0, 0, 0,-2)$&  $x_0^2 x_5 $,  $x_2^3 $ \\
\cline{1-3}
\end{tabular}
\vspace{0.1cm}
\caption{The maximal subsets $M_{\le 0}(\lambda)$}\label{maxsstable}
\end{center}
\end{table}

\begin{remark}\label{remalgo}
Let us briefly indicate the procedure of finding the maximal subsets $M_{\le 0}(\lambda)$. We ask the equivalent question of finding the maximal subsets $M$ of degree $3$ monomials in $5$ variables for which there exists a $1$-PS $\lambda$ such that $M\subseteq M_{\le 0}(\lambda)$. Since the number of subsets $M$ is finite, the computational problem  is finite. However, in this form, it is not effective. We correct this as follows. The normalization assumption on $1$-PS $\lambda$ induces a partial order on the set of monomials (i.e. the order generated by $x_0\ge \dots \ge x_5$, see Mukai \cite[Ch. 7.2, (7.11)]{mukai}).  To determine the maximal subsets $M$ included in some $M_{\le 0}(\lambda)$ it suffices to look at the  maximal monomials (w.r.t. the given partial order) in $M$. Since any two maximal monomials in $M$ are incomparable, it is easy to see that there  are only few possibilities for the set of maximal monomials. This observation gives then  an easy   effective solution. The resulting  subsets  $M=M_{\le 0}(\lambda)$ and the corresponding maximal monomials are given in table \ref{maxsstable} (N.B. the maximal monomials determine both $M$ and $\lambda$).
\end{remark}

\smallskip

Two of the cases from table \ref{maxsstable} are eliminated by the following lemma:
\begin{lemma}\label{elims}
An equation of type S7 or S8 can be reduced to an equation of type S6 by a linear change of coordinates.
\end{lemma}
\begin{proof}
An equation of type S6 has the general form:
$$g_1(x_0,\dots,x_5)=x_0q(x_4,x_5)+f(x_1,\dots,x_5)$$
and it is characterized by the fact that $(1:0:\dots:0)$ is a singular point of corank $3$ for the fourfold defined by $g_1$ (see also \ref{cases6}). To prove the lemma, it suffices to find a singular point of corank at least $3$ for an equation of type S7 or S8. In the case S7, we write the corresponding equation as:
$$g_2(x_0,\dots,x_5)=x_0x_5l(x_1,\dots,x_5)+f(x_1,\dots,x_5).$$
The point $(1:0:\dots:0)$ is a singularity of corank at least $3$ for $g_2$. The reduction of $g_2$ to the form S6 is then obvious.  The case S8 is similar.
\end{proof}

\smallskip

For a cubic form $f$ such that $\mu(f,\lambda)\le 0$, the limit $\lim_{t\to 0} f_t=f_0$ exists and it is invariant with respect to $\lambda$. For the six cases S1--S6, we denote the corresponding invariant parts by $\alpha,\dots,\delta$. A general equation of type $\alpha,\dots,\delta$  has the form:
\begin{eqnarray}
\label{eqalpha}
\ \ \alpha:\  g(x_0,\dots,x_5)&=& x_0\cdot q_1(x_2,\dots,x_5)+x_1\cdot q_2(x_2,\dots,x_5), \\
\label{eqbeta}
\ \ \beta:\ g(x_0,\dots,x_5)&=&ax_0x_4^2+x_0x_5l_1(x_2,x_3)+bx_1^2x_5\\
\notag&&+x_1x_4l_2(x_2,x_3)+f(x_2,x_3),\\
\label{eqgamma}
\ \ \gamma:\ g(x_0,\dots,x_5)&=&x_0q(x_3,x_4,x_5)+x_1^2l_1(x_3,x_4,x_5)\\
\notag&&-2 x_1x_2l_2(x_3,x_4,x_5)+x_2^2l_3(x_3,x_4,x_5),\\
\label{eqdelta}
\delta:\ \ g(x_0,\dots,x_5)&=&x_0\cdot q(x_4, x_5)+f(x_1,x_2,x_3).
\end{eqnarray} 

\begin{remark}
Clearly, $M_0(\lambda)=M_0(\lambda^*)$, where the dual $\lambda^*$ means  opposite weights. Thus, the cases S1--S6 come in pairs. The case S3 is self-dual. The cases S6--S8 can be reduced (and are in fact equivalent) to the single case S6 (cf. lemma \ref{elims}) due to the fact that an equation of type $\delta$ is stabilized  by a $2$-dimensional torus.
\end{remark}

\smallskip

Similarly to case of stable cubics, in order to understand the semistable cubics one has to find the maximal subsets of type $M_{<0}(\lambda)$ as $\lambda$ varies over all normalized $1$-PS. We solve this problem by applying a simple modification of the algorithm described in remark \ref{remalgo}. The results are given in table \ref{maxunstable} below. Each of the maximal subsets $M_{<0}(\lambda)$ is included in some maximal subset $M_{\le 0}(\lambda')$ (see the last column of table \ref{maxunstable}). Again, some of these cases can be removed by using projective equivalences. 

 \begin{lemma}\label{reduceunstable}
An equation of type U7 or U8 can be reduced to an equation of type U1 by a linear change of coordinates. Similarly, an equation of type U9 or U10 can be reduced to an equation of type U5. \qed
\end{lemma}

\begin{table}[htb]
\begin{center}
\renewcommand{\arraystretch}{1.25}
\begin{tabular}[2cm]{|c|l|l|c|}
\hline
Case& Weights of $\lambda$& Maximal monomials& Incl.\\
\hline\hline
(U1)&$(35,23,-1,-13,-19,-25)$&$x_0 x_3 x_5 $, $x_0 x_4^2 $, $x_1 x_2 x_5 $, $x_1 x_3^2 $, $x_2^3 $&(S1)\\
\hline
(U2)&$(47,29,11,-1,-25,-61)$& $x_0 x_2 x_5 $, $x_0 x_4^2 $, $x_1^2 x_5 $, $x_2^2 x_4 $, $x_3^3 $&(S3)\\
\hline
(U3)&$(13, 5, 3,-3,-7,-11)$&$x_0 x_3 x_5 $, $x_0 x_4^2 $, $x_1^2 x_5 $, $x_1 x_3^2 $, $x_2^2 x_4 $&(S4)\\
\hline
(U4)&$(11,11,11,-1,-7,-25)$&$x_0^2 x_5 $, $x_0 x_4^2 $, $x_3^3 $&(S5)\\
\hline
(U5)&$(47,11,-1,-7,-25,-25)$& $x_0 x_4^2 $, $x_1^2 x_4 $,
 $x_1 x_3^2 $, $x_2^3 $&(S6)\\
 \hline
(U6)&$(11,-1,-1,-1,-1,-7)$&$x_0 x_5^2 $, $x_1^3 $&(S6)\\
\hline
(U7)&$(11,11,-1,-7,-7,-7)$& $x_0 x_3^2 $, $x_2^3 $\\
\cline{1-3}
(U8)&$(23,23,-1,-7,-13,-25)$&$x_0 x_2 x_5 $, $x_0 x_4^2 $, $x_2^3 $\\
\cline{1-3}
(U9)&$(29,11,-1,-7,-7,-25)$&$x_0 x_3 x_5 $, $x_1^2 x_5 $, $x_1 x_3^2 $, $x_2^3 $\\
\cline{1-3}
(U10)&$(23,11,-1,-1,-7,-25)$&$x_0 x_2 x_5 $, $x_1^2 x_5 $, $x_1 x_4^2 $, $x_2^3 $\\
\cline{1-3}
\end{tabular}
\vspace{0.1cm}
\caption{Maximal subsets $M_{<0}(\lambda)$}\label{maxunstable}
\end{center}
\end{table}

\smallskip

We introduce the following terminology:
 
\begin{definition}
Let $Y\subset \bP^5$ be a cubic fourfold. We say that {\it $Y$ is of type Sk (Uk)}, if there exists a choice of coordinates on $\bP^5$ such that the defining equation of $Y$ is of type Sk (resp. Uk) as given in table \ref{maxsstable} (resp. table \ref{maxunstable}). A fourfold $Y$ is of type $\alpha$--$\delta$ if its equation can be chosen to be of type (\ref{eqalpha}--\ref{eqdelta}). 
\end{definition}

\smallskip

We conclude the combinatorial aspect of the stability analysis by summarizing the results of the section in the following proposition:
\begin{proposition}\label{sumcomb}
A cubic fourfold is not properly stable (resp. unstable) if and only if it is of type Sk (resp. Uk) for some $k=1,\dots 6$. Furthermore, a strictly semistable cubic fourfold with minimal orbit is of type $\alpha$--$\delta$.
\end{proposition}

\section{Geometric interpretation of stability}\label{sectnonstable}
In this section, we start the geometric analysis of the stability condition for cubic fourfolds. Here, we interpret the results of the tables \ref{maxsstable} and \ref{maxunstable} in terms of the existence of  ``bad flags'' (Mumford \cite[pg. 80]{GIT}).

To explain the occurrence of flags in the geometric analysis, let us recall that a cubic fourfold $Y$ fails to be stable if there exists a $1$-PS $\lambda$ such that $\mu(f,\lambda)\le 0$, where $f$ is the corresponding cubic form. The natural action of the $1$-PS $\lambda$ on the vector space $W=H^0(\bP^5,\calO_{\bP^5}(1))$ determines a weight decomposition $W=\oplus_k W_k$. We define a flag $(F_n)_n$ of subspaces of $W$ by $F_n:=\oplus_{k\le n} W_k$. The condition that $Y$ is not stable with respect to $\lambda$ imposes some specific geometric properties of $Y$ with respect to the associated projective flag $\bP(F_n)$ (e.g. typically $Y$ is singular along the linear subspace corresponding to the first term of the flag). Conversely, certain geometric properties of $Y$ determine a special flag, which in turn leads to existence of a destabilizing $1$-PS. We note that the complexity of the analysis of the cases S1--S6 depends  on the length of the flag  associated to $\lambda$ (i.e. the number of distinct weights of $\lambda$). 

\begin{notation}
Given $\lambda$, we typically denote by $p$, $L$, $P$, $\Pi$, and $H$ the $0,\dots,4$-dimensional members of the associated projective flag.  Under the normalization assumption on $\lambda$, the resulting flag is a  subflag of the standard complete flag on $\bP^5$ (i.e. $p=(1:0:\dots:0), L: (x_2=\dots=x_5=0$, etc.).
\end{notation}

\subsection{Non-Stable Cubic Fourfolds}
The interpretation of the results of table \ref{maxsstable} in terms of destabilizing flags gives  the following geometric characterization of the stable cubic fourfolds:
\begin{theorem}\label{thmsstable}
A cubic fourfold $Y$ is not properly stable iff one of the following conditions holds:
\begin{itemize}
\item[(S1)] $Y$ is singular along a line;
\item[(S2)] $Y$ contains a $3$-plane;
\item[(S3)] $Y$ has a singular point of corank $2$ satisfying the degeneracy conditions of lemma \ref{cases3};
\item[(S4)] $Y$ has a singular point of corank $2$, and it contains the null plane of this singularity;
\item[(S5)] $Y$ is singular along a conic;
\item[(S6)] $Y$ has a singular point of corank $3$ or more.
\end{itemize}
Additionally, a non-stable cubic fourfold $Y$ degenerates to a cubic fourfold of type $\alpha$--$\delta$ (see table \ref{maxsstable}). 
\end{theorem}
\begin{proof} As a consequence of the numerical criterion, we only need to find geometric characterizations for the cubic fourfolds of type S1--S6 (see Prop. \ref{sumcomb}). This is done in lemmas \ref{cases1}--\ref{cases3} below.
\end{proof}

\smallskip

\begin{lemma}[Case S1]\label{cases1}
Let $Y$ be a cubic $4$-fold of type S1. Then $Y$ contains a double line. Conversely, if $Y$ contains a double line then $Y$ is of type S1. 
\end{lemma} 
\begin{proof}
The equation of $Y$ is of type S1, so it can be written as: 
$$g(x_0,\dots,x_5)=x_0q_1(x_2,\dots,x_5)+x_1q_2(x_2,\dots,x_5)+f(x_2,\dots,x_5).$$
The fourfold $Y$ contains the line $L: (x_2=x_3=x_4=x_5=0)$ with multiplicity $2$. The converse is also clear.  Note that the only monomials missing from the above equation  are the $16$ monomials which are not in the square of the ideal $\langle x_2,\dots,x_5\rangle$. 
\end{proof}

\smallskip

The proof of the following lemmas is similar, we omit the details.

\begin{lemma}[Case S2]\label{cases2}
Let $Y$ be a cubic $4$-fold of type S2. Then $Y$ contains a $3$-plane (a three dimensional linear subspace of $\bP^5$). Conversely, if $Y$ contains a $3$-plane then $Y$ is of type S2. Additionally, $Y$ is singular along a curve which is the complete intersection of two quadrics in $\bP^3$. This property is an equivalent geometric characterization of a fourfold of type S2. \qed
\end{lemma} 

\smallskip

\begin{lemma}[Case S4]\label{cases4}
Let $Y$ be a cubic $4$-fold  of type S4. Then $Y$ has a singularity of corank (at least) $2$ and it contains the null plane associated to this singularity. Conversely, if $Y$ contains a singularity of corank $2$ and it also contains the null plane of this singularity, then $Y$ is of type S4.  \qed
\end{lemma} 

\smallskip

\begin{lemma}[Case S5]\label{cases5}
Let $Y$ be a cubic $4$-fold having the defining equation of type S5. Then $Y$ is singular along a (plane) conic. Conversely, if $Y$ is singular along a conic, then $Y$ is of type S5. \qed
\end{lemma} 

\smallskip

\begin{lemma}[Case S6]\label{cases6}
Let $Y$ be a cubic $4$-fold of type S6. Then $Y$ contains a singular point of corank at least $3$. Conversely, if $Y$ contains a singular point of corank $3$ or more then  $Y$ is of type S6. \qed
\end{lemma} 

\smallskip

In all the cases except S3, the destabilizing $1$-PS $\lambda$ has at most $3$ different weights. Thus, the flag  associated  to $\lambda$ is rather simple. As a consequence,  the analysis of the previous lemmas is  straightforward. In contrast, in the case S3, the associated $1$-PS has $5$ distinct weights, producing a flag $p\in L\subset \Pi\subset H\subset \bP^5$, where $p$ is point, $L$ a line, $\Pi$ a $3$-plane, and  $H$ a hyperplane. It follows that the geometry of a cubic fourfold of type S3 is quite complicated, reflecting the complexity of the flag.

\begin{lemma}[Case S3]\label{cases3}
Let $Y$ be a cubic $4$-fold  of type S3. Assume that $Y$ is not of type S6. Then $Y$ has a singularity $p\in Y$ of corank $2$ with the property that there exists a hyperplane $H$ such that:
\begin{itemize}
\item[i)]  the intersection of $H$ with the projective tangent cone  to  $p$ (a rank $3$ quadric in $\bP^5$)
 is a $3$-plane $\Pi$ (with multiplicity $2$); 
\item[ii)]  the restriction of $Y$ to $H$ contains a double line $L$ passing through $p$;
\item[iii)] the restriction of $Y$ to $\Pi$ consists of $3$ planes meeting in the line $L$.  
\end{itemize}
Conversely, if a cubic fourfold $Y$ has these properties then $Y$ is of type S3. 
\end{lemma} 
\begin{proof}
A fourfold $Y$ of type S3 has the defining equation:
$$g(x_0,\dots,x_5)=x_0Q(x_1,\dots,x_5)+F(x_1,\dots,x_5)$$
where $Q(x_1,\dots,x_5)=ax_4^2+x_5l(x_2,x_3,x_4,x_5)$ and
$$F(x_1,\dots,x_5)=bx_1^2x_5+x_1(x_2l_1(x_4,x_5)+x_3l_2(x_4,x_5)+q(x_4,x_5))+f(x_2,\dots,x_5).$$
The quadric $Q$ has rank at most $3$. It follows that $p=(1:0:\dots:0)\in\bP^5$ is a singular point of $Y$ of corank $2$ or more. We assume that $Y$ is not of type S6. In particular, the corank at $p$ is exactly $2$ (i.e. $a\neq 0$, and $l(x_2,x_3,x_4,x_5)$, $x_4$ and $x_5$ are linearly independent). The null plane of the singularity at $p$ is given by $P: (l(x_2,x_3,x_4,x_5)=x_4=x_5=0)$. 

The $1$-PS $\lambda$ associated to $Y$ of type S3 singles out the line  $L: (x_2=x_3=x_4=x_5=0)$,
the $3$-plane $\Pi:(x_4=x_5=0)$, and the hyperplane $H:(x_5=0)$. Together with $p$ and $P$ defined above, we obtain a 
 full flag $p\in L\subset P\subset \Pi\subset H$. The following  properties of the flag are immediate:
 \begin{itemize}
 \item[i)] $H$ meets the projective tangent cone to $p$ (given by the quadric $Q$ in $\bP^5$) in the $3$-plane $\Pi$;
 \item[ii)] the line $L$ is a double line for the cubic $3$-fold $X=Y\cap H\subset H\cong \bP^4$ (note that all monomials of $g$ are in the ideal $\langle x_2,x_3,x_4\rangle^2+\langle x_5\rangle$); 
 \item[iii)] the restriction of $Y$ to $\Pi\cong \bP^3$ consists of three planes meeting in the line $L$  (note that  the restriction is given by $(h(x_2,x_3)=0)\subset \Pi\cong \bP^3$, where $h(x_2,x_3)=f(x_2,x_3,0,0)$). 
\end{itemize}
We obtain one direction of the lemma. 
  
 Conversely, we assume that $Y$ has the geometric properties stated in the lemma. These properties  recover the flag $p\in L\subset \Pi\subset H$, which in turn  allows us to find coordinates such that $Y$ has an equation of type S3. Specifically,  let $p=(1:0:\dots:0)\in Y$ be the corresponding singular point. Since $p$ has corank $2$, we can write the equation of $Y$ as:
$$g(x_0,\dots,x_5)=x_0(x_4^2+x_5x_3)+F(x_1,x_2,x_3,x_4,x_5).$$
The null plane of the singularity at $p$ is given by $P: (x_3=x_4=x_5=0)$. Let $H$ be a hyperplane as in the lemma. The assumption i)  on $Y$ forces $H$ to contain the null plane $P$. Furthermore, we can view $H$ as being spanned by $P$ and a tangent line to the conic  $(x_4^2+x_3x_5=0)$ in the plane $P':(x_0=x_1=x_2=0)$. In particular,  we get a partial flag $p\in  P\subset \Pi\subset H$, where $\Pi$ is the reduced intersection of $H$ with the projective tangent cone at $p$. By a change of variables involving only $x_3$,  $x_4$, and $x_5$, and fixing the conic $x_4^2+x_3x_5$  we can assume that $\Pi$ is given by $\Pi:(x_4=x_5=0)$ and $H$ is given by $H:(x_5=0)$. The assumption ii) on $Y$, implies in particular that the line $L$ through $p$ lies in the null plane $P$. Thus, we obtain a full flag $p\in L\subset  P\subset \Pi\subset H$. By a change of variables involving only $x_1$ and $x_2$ we can further assume that $L$ is given by $L:(x_2=x_3=x_4=x_5=0)$. With these normalizations, it is easily seen that the condition iii) is equivalent to vanishing of the coefficients of the monomials: $x_1^3$, $x_1^2x_2$, $x_1^2x_3$, $x_1x_2^2$, $x_1x_2x_3$ and $x_1x_3^2$ in the cubic $F(x_1,\dots,x_5)$.
Thus, we write the cubic $F$ as 
\begin{eqnarray*}
F(x_1,\dots,x_5)&=&x_1^2(bx_5+cx_4)+x_1(x_2l_1(x_4,x_5)+x_3l_2(x_4,x_5)+q(x_4,x_5))\\
                          &&+f(x_2,x_3,x_4,x_5).
\end{eqnarray*}
Using again the condition ii), we obtain that the coefficient $c$ in the equation of $F$ vanishes. In conclusion, $Y$ is of type S3.
\end{proof}

\subsection{Unstable cubic fourfolds}\label{sectunstable}
Theorem \ref{thmsstable} gives a satisfactory geometric description of the stable locus for cubic fourfolds. In principle, a similar analysis  of the results of table \ref{maxunstable} would give  a separation of the strictly semistable cubics from the unstable ones (e.g. see the results of Allcock \cite[Thm. 1.3, 1.4]{allcock1} for cubic threefolds). Unfortunately, we find the computations  tedious, without a clear pattern. Thus, we choose to give only two results that are relevant for the subsequent discussion.

\begin{lemma}[Case U6]\label{caseu6}
Let $Y$ be a cubic $4$-fold of type U6. Then $Y$ contains a singular point of corank at least $4$. Conversely, if $Y$ contains a singular point of corank $4$ or more then  $Y$ is of type U6. In particular a cone over a cubic threefold is unstable.
\end{lemma}
\begin{proof}
This is the same as lemmas \ref{elims} and \ref{cases6}.
\end{proof}

\smallskip

\begin{lemma}
Assume that $Y$ is a fourfold having  only   singularities of type $A_n$ or $D_n$ (including $A_{\infty}$ and $D_{\infty}$). Then $Y$ is semi-stable.
\end{lemma}
\begin{proof}
If $Y$ is unstable, it is of one of the types U1--U6. In all these cases, except U4, there exists a singularity at $p=(1:0:\dots:0)$ of corank at least $3$, or of corank $2$ and third jet $x_1^3$. It follows that the singularity at $p$ is worse than $D_{\infty}$. In the case U4,  $Y$ is singular along a conic with singularities of $A_2$ transversal type (i.e. the singularities of a generic hyperplane section of $Y$ are of type $A_2$). Thus, the singularities of $Y$ are worse than $D_\infty$.
\end{proof}

\smallskip

\begin{remark}
For a fixed cubic fourfold $Y$ of equation $f$, it is not hard to decide if it is unstable or not. Namely,  $Y$ is unstable if and only if it satisfies one of the geometric conditions of theorem \ref{thmsstable} and additionally the associated limit $f_0=\lim_{t\to 0} \lambda(t)f$ (where $\lambda$ is determined by Thm. \ref{thmsstable}) is unstable (N.B. $f_0$ is of type $\alpha$--$\delta$ and as such it is  covered by lemmas \ref{casealpha}--\ref{casedelta}). 
\end{remark}

\section{The minimal orbits and their normal forms}\label{sectminorbits} 
The geometric invariant theory compactifies the geometric quotient $\calM^s$ by adding boundary components that parametrize strictly semistable cubic fourfolds with minimal orbits.  By proposition \ref{sumcomb}, we know that such a cubic fourfold is of type $\alpha$--$\delta$. It follows that there are four boundary components for $\overline{\calM}\setminus \calM^s$, which, by abuse of notation, we label also by $\alpha$--$\delta$. In this section, we discuss which of the fourfolds of type $\alpha$--$\delta$ are in fact semistable, and  the structure of the resulting boundary components. As noted in \S\ref{sectunstable}, a direct approach based on the numerical criterion is cumbersome. We use instead the following criterion of Luna \cite[Cor. 1]{luna}: {\it Let  $X$ be an affine $G$-variety, and $x\in X$ a point stabilized by a reductive subgroup $H$. Then the orbit $G\cdot x$ is closed if and only if the orbit $N_G(H)\cdot x$ is closed in $X^H$.}

\smallskip

 We apply Luna's criterion for the affine space $X=\Sym^3(W)$ and the connected component $H$ of the stabilizer of a general cubic fourfold of type  $\alpha$--$\delta$.  The advantage of this approach is that it allows us to  check that an orbit of type $\alpha$--$\delta$ is closed by  using smaller subgroups of $G=\SL(W)\cong \SL(6)$.  In other words, in order to understand the boundary it suffices to study the simpler action of $N^G(H)/H$ on $X^H$ instead of the full action of $G$ on $X$ (N.B.  $H$ acts trivially on $X^H$). We note also that the natural morphism $X^H/N_G(H)\to X/G$ is finite (Luna \cite[Main Thm.]{luna}) and that Luna's criterion can be applied by using the centralizer $C_G(H)$ instead of the normalizer $N_G(H)$ (see \cite[pg. 221-222]{popov}). As an immediate application of Luna's results, we obtain the dimensions of the boundary  strata.    
 
\begin{lemma}
A generic cubic fourfold of type $\alpha$--$\delta$ gives a closed orbit. Therefore, each of the types $\alpha$--$\delta$ gives an irreducible boundary component for the GIT compactification $\calM^s\subset \overline{\calM}$. The dimensions of these boundary strata are $1$ for $\alpha$ and $\delta$, $2$ for $\gamma$, and $3$ for $\beta$. \qed
\end{lemma}
 
 \smallskip
 
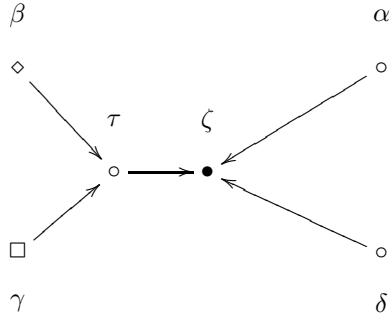
\begin{figure}[htb]
$$\xymatrix@R=.25cm{
{\beta}   &&                      &&{\alpha}             \\
{\diamond}\ar@{->}[ddr]&&                      &&{\circ}\ar@{->}[ddll]                \\
 &{\tau}                        &{\zeta}\\
&{\circ}\ar@{->}[r]   &{\bullet}        \\
  &&                      &&                    \\
{\Box}\ar@{->}[uur]&&                      &&{\circ}\ar@{->}[uull]                  \\
{\gamma}   &&                      &&{\delta}                     \\
}
$$
\caption{Incidence of the  boundary components of $\calM^s$ in $\overline{\calM}$}\label{gitboundary}
\end{figure}

Before starting the detailed analysis of the strata $\alpha$--$\delta$, we observe that the  types $\beta$ and $\gamma$ have as a common specialization the curve:
\begin{eqnarray}\label{eqtau}
\ \ \ \ \ \ \ \tau:\  \ g(x_0,\dots,x_5)&=& -\left| \begin{array}{ccc}
x_0&x_1&a\cdot x_2\\
x_1&x_5&x_3\\
x_2&x_3&x_4
\end{array}\right|\\
\notag&=&x_0(x_3^2-x_4x_5)+x_1^2x_4- (a+1) x_1x_2x_3+a x_2^2x_5,\end{eqnarray}
where $a\in \bC$. The stabilizer of an equation of type $\tau$ contains the $2$-dimensional torus generated by $\lambda: (2,1,0,-1,-2,0)$ and $\lambda': (4,1,1,-2,-2,-2)$. The equation (\ref{eqtau}) further degenerates (for $a=0$) to:
\begin{equation}\label{eqzeta}
\zeta:\ \ g(x_0,\dots,x_5)=
\left| \begin{array}{ccc}
x_0&x_1&0\\
0&x_5&x_3\\
x_2&0&x_4
\end{array}\right|=x_0x_4x_5+x_1x_2x_3,
\end{equation}
which is stabilized by a $4$-dimensional torus. The orbit $\zeta$ is also a specialization of the cases $\alpha$ and $\delta$. The resulting incidence diagram is given in figure \ref{gitboundary}. 

\begin{remark}\label{defomega} 
A general fourfold of type $\tau$ is singular along $3$ conics meeting pairwise in a point. There exists a special point $\omega$ on $\tau$, corresponding to $a=1$ in (\ref{eqtau}), which gives the determinantal cubic fourfold:
\begin{eqnarray}\label{eqomega}
\omega:\ \ g(x_0,\dots,x_5)&=&-\left| \begin{array}{ccc}
x_0&x_1&x_2\\
x_1&x_5&x_3\\
x_2&x_3&x_4
\end{array}\right|\\
\notag &=&x_0(x_3^2-x_4x_5)+x_1^2x_4-2x_1x_2x_3+x_2^2x_5.
\end{eqnarray}
The determinantal cubic is singular along the Veronese surface and is stabilized by a subgroup of $\SL(6)$ isomorphic to $\SL(3)$. 
\end{remark}

\smallskip

Our first result regarding the minimal orbits establishes the semi-stablity of the degenerate cases $\tau$ and $\zeta$. 
\begin{lemma}\label{casetau}
A fourfold of type $\zeta$ is semi-stable with closed orbit. A fourfold $Y$ of type $\tau$ (in particular $\omega$) is semi-stable. If $a\neq 0$, the orbit of $Y$ is closed. For $a=0$, the closure of the orbit of $Y$ contains the orbit of $\zeta$.
\end{lemma}
\begin{proof}
This follows from Luna's criterion cited above. The stabilizer of (\ref{eqtau}) contains a $1$-PS $H$ of distinct weights (e.g. $(6,2,1,-3,-4,-2)=\lambda\cdot \lambda'$). Thus it suffices to check the semi-stability with respect to the standard maximal torus $T=C_G(H)$ in $G$. The proposition follows easily. For example, the fact that $\zeta$ is semi-stable is equivalent  to saying that $a_0+\dots+a_5=0$ implies that either $a_0+a_4+a_5\ge 0$ or $a_1+a_2+a_3\ge 0$, where $(a_0,\dots,a_5)$ are the weights of a $1$-PS of $T$.
\end{proof}

\smallskip

We now do the case-by-case analysis of the minimal orbits of type $\alpha$--$\delta$. The common feature of all these cases is that the analysis reduces to some well-known  lower dimensional GIT problem. For example, the case $\alpha$ reduces to the analysis of the stability for pencils of quadrics (see  \cite{wallquadrics} and \cite{avritzermiranda}).

\begin{lemma}\label{casealpha}
Let $Y$ be a fourfold of type $\alpha$ given by the equation (\ref{eqalpha}). Regard the equation of $Y$ as a pencil (parametrized by the double line on $Y$)  of quadrics in the variables $x_2,\dots,x_5$. Let $\Delta(x_1,x_2)$ be the discriminant (a binary quartic) of the pencil, and $E$ the base locus of the pencil. One of the following holds:
\begin{itemize}
\item[i)] If the roots of $\Delta(x_1,x_2)$ are distinct, then $Y$ gives a minimal orbit and the equation of $Y$ can be taken as:  
\begin{equation}\label{eqalpha2}
g(x_0,\dots,x_5)= x_0\cdot (a x_3^2+x_4^2+x_5^2)-x_1\cdot (x_2^2+bx_3^2+x_4^2)
\end{equation}  
with $\lambda=\frac{a}{b}\neq 0,1,\infty$. The fourfold $Y$  is singular along $E$  and along the line $L: (x_2=\dots=x_5=0)$. The singularities along $E$ are of type $A_{\infty}$. Similarly, those along $L$ are of type $A_\infty$ at all but $4$ points (the zeros of $\Delta$), where they become of type $D_{\infty}$. 
\item[ii)] If $\Delta(x_1,x_2)$ has a double root, then  $Y$ degenerates to $\zeta$.
\item[iii)] If $\Delta(x_1,x_2)$ has a triple root (or vanishes), then $Y$ is unstable. 
\end{itemize}
\end{lemma}
\begin{proof}
The connected component of the stabilizer of an equation of type $\alpha$ is the $1$-PS $H$ of weights $(2,2,-1,-1,-1,-1)$. It follows that $C_G(H)$ acts on an equation of type $\alpha$ as $\SL(2)\times \SL(4)$:  $\SL(2)$ acts on the variables $x_0$ and $x_1$, and $\SL(4)$ acts on $x_2,\dots,x_5$. More intrinsically, the action of $C_G(H)/H$ (N.B. $H$ is abelian) on $X^H$ is equivalent to the natural action of  $\SL(2)\times \SL(4)$ on $V\otimes \Sym^2(U)$, where $V$ and $U$ are the standard representations of $\SL(2)$ and $\SL(4)$ respectively. Therefore, the GIT analysis for a fourfold of type $\alpha$ is equivalent to the GIT analysis for the pencil of quadric surfaces in $\bP^3$ given by $q_1$ and $q_2$ (see \cite[\S2]{wallquadrics}). According to Wall \cite[\S4(a)]{wallquadrics}, the pencil is semi-stable if and only if the multiplicity of the roots of $\Delta(x_1,x_2)$ is at most $2$. If the roots of $\Delta$ are distinct, the corresponding orbit is closed and 
the quadrics $q_1$ and $q_2$ can be simultaneously diagonalized, giving the equation  (\ref{eqalpha2}) (see \cite[Prop. 2]{avritzermiranda}). The invariant of the pencil is the base locus, the elliptic curve $E$.   If $\Delta$ has a double root, then the corresponding orbit contains in its closure the orbit $\zeta$.
\end{proof}

\smallskip

The case $\delta$ is quite similar, we omit the details.
\begin{lemma}\label{casedelta}
Let $Y$ be a fourfold of type $\delta$  (see (\ref{eqdelta})). One of the following holds:
\begin{itemize}
\item[i)] if $q(x_4,x_5)$ has a double root, then $Y$ is unstable of type U6;
\item[ii)] if $f(x_1,x_2,x_3)$ has a cusp, then $Y$ is unstable of type U5;
\item[iii)] if $f(x_1,x_2,x_3)$ has a node, then  $Y$ degenerates to $\zeta$;
\item[iv)] if none of the items i-iii are true, then $Y$ is semistable with closed orbit. After a linear change of coordinates we can assume $Y$ is given by:
\begin{equation}\label{eqdelta2}
g(x_0,\dots,x_5)=x_0x_4x_5+f(x_1,x_2,x_3)
\end{equation}
where $f(x_1,x_2,x_3)$ is a smooth plane cubic.
The fourfold $Y$ has $3$ singularities of type $\widetilde{E}_6$ (of the same modulus and non-colinear).  \qed
\end{itemize}
\end{lemma}

\smallskip

As the dimension of stratum increases, the analysis gets more involved.
\begin{lemma}\label{casegamma}
Let $Y$ be a fourfold of type $\gamma$ given by (\ref{eqgamma}). Assume that the rank of the quadric $q$ is $3$ (otherwise is of type S6 and can be treated as such).  Let $q_2(x_3,x_4,x_5)=\left| \begin{array}{cc} l_1&l_2\\ l_2&l_3\end{array}\right|$. Then one of the following holds:
\begin{itemize}
\item[i)] The conics $q$ and $q_2$ meet in $4$ distinct points. Then $Y$ gives a minimal orbit and its equation  can be taken as 
\begin{equation}\label{eqgamma2}
g(x_0,\dots,x_5)= x_0(x_3^2-x_4x_5)+x_1^2x_4-2  x_1x_2 l(x_3,x_4,x_5) +x_2^2x_5
\end{equation} 
with the conics $q=(x_3^2-x_4x_5)$ and $q_2=x_4 x_5-l(x_3,x_4,x_5)^2$ meeting in four distinct points.
The fourfold $Y$ is singular in the point $(1:0:\dots:0)$ with a singularity of type $\widetilde{E_7}$ and along the conic given by $q$ in the $2$-plane $(x_0=x_1=x_2=0)$.  The singularities along the conic are $A_\infty$, except at the $4$ special points, where they are $D_{\infty}$. 
\item[ii)] If the conics $q$ and $q_2$ do not meet transversally, then the orbit of $Y$ is not closed. Either $Y$ is semistable and degenerates to a fourfold of type $\tau$ or  $Y$ is unstable. 
\end{itemize}
\end{lemma}
\begin{proof}
The connected component of the stabilizer of a fourfold of type $Y$ is the $1$-PS $H$ of weights $(4,1,1,-2,-,2,-2)$. Up to the scaling of the variables, $C_G(H)$ acts on the stratum $\gamma$ as $\SL(2)\times \SL(3)$ with $\SL(2)$ acting on $x_1$ and $x_2$ and $\SL(3)$ on $x_3,\dots, x_5$. The action of  $\SL(2)$ changes the matrix $\left(\begin{array}{cc} l_1&l_2\\ l_2&l_3\end{array}\right)$ to a conjugate (N.B. $\SL(2)$ acts trivially on $q$). Therefore, the invariant part is the determinant $q_2$. We obtain two conics $q$ and $q_2$ in $\bP^2$. By acting with $\SL(3)$ we can bring the two conics to the normal form (\ref{eqgamma2}). The remaining part follows by an analysis of the possible degenerations. 
\end{proof}

\smallskip

The last remaining case is the $3$-dimensional boundary stratum $\beta$.
\begin{lemma}\label{casebeta}
Let $Y$ be a fourfold of type $\beta$ given by (\ref{eqbeta}). The orbit of $Y$ is not closed iff any of the following conditions are satisfied: $a$, $b$, $l_1$, or $f$ vanish, $l_1$ is a factor of $f$, or $l_2^2$ is a factor of $f$ (including the case $f$ has a double root and $l_2$ vanishes).  If the orbit of $Y$ is not closed, then either $Y$ is unstable or $Y$ degenerates to a fourfold of type $\tau$ (including $\zeta$). If $Y$ is semistable with closed orbit, then: 
\begin{itemize}
\item[i)] Generically,  $Y$ has two isolated singularities of type $\widetilde{E}_8$ (of same modulus).
\item[ii)] If $Y$ doesn't have isolated singularities, then $Y$ is singular along a rational normal curve of degree $4$ (the case labeled $\sigma$ (\ref{eqsigma}) and discussed in \S\ref{nonisolated2}). 
\end{itemize}
\end{lemma}
\begin{proof}
Let $H$ be the $1$-PS $(2,1,0,0,-1,-2)$ that stabilizes $Y$. The center of $C_G(H)$ is a $4$-torus acting on an equation of type $\alpha$ by scaling of the variables. The factor group $C_G(H)/Z(C_G(H))$ is isomorphic to $\PGL(2)$ acting  on  $l_1$, $l_2$, and $f$ by linear changes of the variables $x_2$ and $x_3$. This makes the structure of $X^H/C_G(H)$ (and consequently that of strata $\beta$) quite transparent. In particular, in order to understand the closed orbits, we can apply a slight modification of the computer program used in \S\ref{sectgit}. As a result, we obtain the degeneracy conditions of the lemma (N.B. most of these conditions can be seen directly).
 
 By inspecting  (\ref{eqbeta}), it is easy to see that  there are two singularities of corank $2$ at $(1:0:\dots:0)$ and $(0:\dots:0:1)$. In the generic case, they are of type $\widetilde{E_8}$ (see also section \ref{sectisolated}, esp. the proof of \ref{deforme7e8}).  Due to the obvious symmetry, the singularities have the same modulus. If $Y$ does not have isolated singularities, then $Y$ is singular along a curve $C$. We distinguish two cases, either $C$ is contained in a linear $3$-dimensional subspace of $\bP^5$ or not. The former case gives a non-closed orbit which degenerates to the case $\tau$ (see also \S\ref{nonisolated1}). Thus, by theorem \ref{thmsingularities}, we can assume that  $C$ is either a rational normal curve of degree $4$ or an elliptic normal curve of degree $6$. Since $Y$ and (consequently) $C$ are stabilized by a subgroup of $\SL(6)$ isomorphic to $\bC^*$, the elliptic case is excluded. We conclude that $C$ is a rational normal curve of degree $4$. \end{proof}

\section{Stability of cubic fourfolds with isolated singularities}\label{sectisolated}
In this and following section, we refine the results of section \ref{sectnonstable}, by interpreting the failure of stability for cubic fourfolds in terms of the existence of bad singularities (rather than bad flags). The main tool, as noticed  by  Allcock \cite[\S2]{allcock1}, is the following result (see  \cite[pg. 209]{AGV1}) that  identifies the analytic type  of certain classes of singularities: 

\begin{theorem}[Arnold]\label{recognsing}
If an analytic function $f(x_1,\dots,x_n)$ is semiquasihomogeneous with respect to the weights corresponding to an $A_n$, $D_n$, $E_r$ or  $\widetilde{E_r}$ (for $r=6,7,8$) singularity, then $f$ has a singularity of that type at the origin.
\end{theorem}

We recall that  a convergent power series $f(x_1,\dots,x_n)$ is {\it semiquasihomogeneous} with respect to a choice of weights if the leading term $f_0(x_1,\dots,x_n)$ (w.r.t. the weighting) defines an isolated singularity at the origin. The weights that make the A-D-E singularities quasihomogeneous are well known (e.g.  \cite[\S2]{allcock1}). Similarly, those for $\widetilde{E}_r$ are $(\frac{1}{3},\frac{1}{3},\frac{1}{3},\frac{1}{2},\dots,\frac{1}{2})$,  $(\frac{1}{4},\frac{1}{4},\frac{1}{2},\dots,\frac{1}{2})$, $(\frac{1}{3},\frac{1}{6},\frac{1}{2},\dots,\frac{1}{2})$ for $r=6,7,8$ respectively. 

\smallskip

An immediate consequence of the previous theorem is the following lemma:

\begin{lemma}\label{genericer}
A general cubic fourfold $Y$ of type Sk, where $k=3$, $4$, or $6$, has a singularity of type $\widetilde{E_r}$, where $r=8$, $7$, or $6$ respectively. \qed
\end{lemma}

It follows that the cubic fourfolds with very mild singularities are stable.
\begin{corollary}\label{corsimplestable}
A cubic fourfold having at worst simple  singularities is stable.
\end{corollary}
\begin{proof}
Let $Y_0$ be such a fourfold. Assume that $Y_0$ is not properly stable. According to proposition \ref{sumcomb}, the fourfold $Y_0$ is of one of the types S1--S6. If the type is S1, S2, or S5, then $Y_0$ is singular along a curve (cf. theorem \ref{thmsstable}), thus a contradiction. It remains to consider the cases Sk for k=3, 4, or 6. Choose coordinates such  that the equation $g$ of $Y_0$ is as given by table \ref{maxsstable}. It follows that $p=(1:0\dots,0)\in Y_0$ is a singular point of $Y_0$. By considering a general deformation of $Y_0$ that preserves the Sk type (w.r.t. the given choice of coordinates), we see that the singularity at $p$ deforms to a singularity of type  $\widetilde{E}_r$ (cf. \ref{genericer} above).  This is a contradiction to the well known fact that no A-D-E singularity deforms to an $\widetilde{E}_r$ singularity.  We conclude that a cubic fourfold simple (isolated) singularity is stable.\end{proof}

The following lemmas establish a converse to the previous corollary:
\begin{lemma}\label{deforme6}
Let $Y$ be a cubic fourfold having a singularity that deforms to $\widetilde{E}_6$. Then $Y$ is of type S6. In particular, $Y$ is not stable.  
\end{lemma}
\begin{proof}
A singularity that deforms to $\widetilde{E}_6$ has corank at least $3$ (by the semi-continuity of the corank). Therefore, $Y$ is of type S6 (cf. lemma \ref{cases6}). 
\end{proof}

\begin{lemma}\label{deforme7e8}
Assume that the cubic fourfold $Y$ has an isolated singularity of corank $2$ that deforms to $\widetilde{E}_8$ or $\widetilde{E}_7$, then $Y$ is of type S3 or S4. In particular, $Y$ is not stable. 
\end{lemma}
\begin{proof}
An equivalent formulation of the lemma is:   if $p\in Y$ is a singularity of corank $2$, then either $p$ is of type $D_n$ or $E_r$, or   $Y$ is of type S3 or S4.  This follows from a careful analysis of the position of null plane of the singularity at $p$  and a systematic application of theorem \ref{recognsing}, as sketched below. 

Without loss of generality, we assume that the  singularity $p\in Y$ of corank $2$ is at $(1:0\dots:0)\in \bP^5$. For an appropriate choice of coordinates we can write the equation $g$ of $Y$ as: $g(x_0,\dots,x_5)=x_0(x_4^2-x_3x_5)+F(x_1,\dots,x_5)$ (N.B. the rank of the tangent cone is $3$).  We project $Y$ onto a hyperplane, say $(x_0=0)$, not passing through the singularity. In this hyperplane we take the coordinates $(x_1:\dots:x_5)$, and denote by $Q$ the projectivized tangent cone  (the quadric given by $x_4^2-x_3x_5$), and by $X$ the cubic threefold given by $F(x_1,\dots,x_5)$. The null plane $P$ of the singularity projects to the line $\bar{P}:(x_3=x_4=x_5=0)$, the singular locus of the quadric $Q$. We have the following possibilities for the position of the line $\bar{P}$ relative to $X$:
\begin{itemize}
\item[(a)] $\bar{P}$ is transversal to $X$;
\item[(b)] $\bar{P}$ meets $X$ with multiplicity $2$ in a point;
\item[(c)] $\bar{P}$ meets $X$ with multiplicity $3$ in a point;
\item[(d)] $\bar{P}$ is contained in $X$.
\end{itemize}

In the first two cases,  $Y$ has the third jet either $x_1^3+x_2^3$ or $x_1x_2^2$. Thus, the singularity at $p$ is   of type $D_n$ (see \cite[pg. 190]{AGV1}). For instance, in the case (a) the equation (in affine coordinates) of the singularity at $p$ can be taken to be: 
$$x_4^2-x_3x_5+x_1^3+x_2^3+(\textrm{higher order terms}).$$ 
Thus, $p$ is of type $D_4$ (cf. theorem \ref{recognsing}). 

The case (d) is equivalent to saying that the null plane $P$ of the singularity at $p$ is contained in the fourfold $Y$. By lemma \ref{cases4},  $Y$ is of type S4 and we are done. 

The only remaining case is (c). In this situation, either $Y$ has a singularity of type $E_r$ for $r=6,7,8$ at $p$, or it satisfies the degeneracy conditions of lemma \ref{cases3}.  First,  one verifies that unless both conditions i) and ii) of lemma \ref{cases3} (they impose a single additional condition to (c)) are satisfied then $Y$ has a singularity of type $E_6$ at $p$. Given that i) and ii) are satisfied, the condition iii) differentiates between the  case when singularity at $p$ is $E_7$ or $E_8$ and the case when the singularity is $\widetilde{E}_8$ (or worse).  Specifically, the projection of the $3$-plane $\Pi$ of the lemma \ref{cases3} is a $2$-plane $\bar{\Pi}$ which cuts on $X$ a singular plane cubic curve $C\subset \bar{\Pi}$. If $C$ is nodal or cuspidal the singularity at $p$ is of type $E_7$ or $E_8$ respectively (again by applying theorem \ref{recognsing}).  Due to some geometric restrictions on $C$ imposed by i) and ii), the only possibility for $C$ to be reducible is  to consist of $3$ lines meeting in one point. This situation is equivalent to the condition iii) of the lemma \ref{cases3}, and generically (if the three lines are distinct) gives an $\widetilde{E}_8$ singularity. We obtain that, in the case (c), either the singularity at $p$ is of type $E_r$ or $Y$ is of type S4, concluding the proof of the lemma.
\end{proof}

We conclude:
\begin{theorem}\label{thmsimplesing}
Let $Y$ be a cubic fourfold with isolated singularities. Then $Y$ is stable if and only if  it has at worst simple singularities. 
\end{theorem}
\begin{proof}
One direction is given by \ref{corsimplestable} above. Conversely, any  non-simple isolated singularity deforms to a singularity of type $\widetilde{E}_r$ singularities. Thus, the result follows from  lemmas \ref{deforme6} and \ref{deforme7e8}.
\end{proof}

\begin{remark} We make the following remarks on the extent to which the stability can be characterized in terms of singularities.
\begin{itemize} 
\item[i)] An important reason why a result such as theorem \ref{thmsimplesing} is possible is that much of the information needed to decide that the singularity is simple is given by corank and third jet. On the other hand, in the case of cubics, the same invariants are the first that come up in the GIT analysis.
\item[ii)] It is not hard to extend the results of this section to include also the singularities $\widetilde{E}_r$. We can show that a cubic fourfold with at worst $\widetilde{E}_r$ is semi-stable. Additionally, if $Y$ contains a singularity of type $\widetilde{E}_{r}$, then  $Y$ degenerates to a fourfold of type $\beta$, $\gamma$, or $\delta$ for $r=8$, $7$, or $6$ respectively. 
\item[iii)] On the negative side, it does not seem to us that it is possible to extend the results beyond the   $\widetilde{E}_r$ singularities. In particular, it is unlikely that reasonable analogous statements to those of Allcock \cite{allcock1} (on the classification of unstable cubics in terms of singularities, esp. \cite[Thm. 1.3 (iii)]{allcock1} and \cite[Thm. 1.4 (iv)]{allcock1}) are possible for cubic fourfolds.
\end{itemize}
\end{remark}

\section{Cubic fourfolds with non-isolated singularities}\label{sectnonisolated}
In order to complete the geometric analysis of stability for cubic fourfolds, we have to understand the cubic fourfolds with non-isolated singularities. More precisely, we are interested in finding the cubic fourfolds with non-isolated singularities that are stable, and a parameterization for them. This is answered by the following theorem that summarizes the results of the section. 

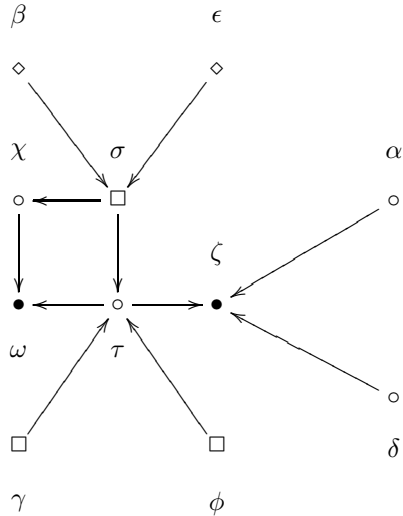
\begin{figure}[htb]
$$\xymatrix@R=.25cm{
&{\beta}   && {\epsilon}                     &&       \\
&{\diamond}\ar@{->}[dddr]&&   {\diamond}\ar@{->}[dddl]       &&          \\
\\
&{\chi}&{\sigma}   &&                      &{\alpha}             \\
&{\circ}\ar@{->}[dd]  \ar@{<-}[r]&{\Box}\ar@{->}[dd]&                      &&{\circ}\ar@{->}[ddll]                \\
 &&                      &{\zeta}\\
&{\bullet}\ar@{<-}[r]&{\circ}\ar@{->}[r]                  &{\bullet}        \\
&{\omega}  &{\tau}  &                      &&                    \\
& && && {\circ}\ar@{->}[uull] \\
&{\Box}\ar@{->}[uuur]&&{\Box}\ar@{->}[uuul]                      && {\delta}                     \\
&{\gamma}   &&{\phi}                      &&                 \\
}
$$
\caption{The incidence of boundary components of $\calM$ in $\overline{\calM}$}\label{gitboundary2}
\end{figure}

\begin{theorem}\label{thmsingularities}
Let $Y$ be a cubic fourfold with non-isolated singularities. Then one of the following holds:
\begin{itemize}
\item[i)] The singular locus of $Y$ contains a line, in which case the orbit of $Y$ contains in its closure an orbit of type $\alpha$.
\item[ii)] The singular locus of $Y$ contains  a conic, in which case the orbit of $Y$ contains in its closure an orbit of type $\gamma$.
\item[iii)] The singular locus of $Y$ contains an elliptic normal quartic curve or a degeneration of it, in which case the orbit of $Y$ contains in its closure an orbit of type $\alpha$.
\item[iv)] The singular locus of $Y$ contains a rational normal quartic curve. A general fourfold of this type  is stable. The closure in $\overline{\calM}$ of the locus  of stable cubic fourfolds singular along a rational normal quartic is an irreducible threefold $\epsilon$. 
\item[v)] The singular locus of $Y$ contains an elliptic normal curve of degree $6$ or a degeneration of it.  A general fourfold of this type  is stable.  The closure in $\overline{\calM}$ of the locus  of stable cubic fourfolds singular along a rational normal sextic is an irreducible surface $\phi$. 
\end{itemize}
Furthermore, the threefold $\epsilon$ meets the GIT boundary $\overline{\calM}\setminus \calM^s$ along a surface $\sigma$, which is contained in the stratum $\beta$. The surface $\phi$ meets the boundary along the curve $\tau$. The adjacency diagram of the strata $\alpha,\dots,\phi$ is given in figure \ref{gitboundary2}.
\end{theorem}
\begin{proof}
Let $C$ be an irreducible curve included in the singular locus  of $Y$. If $C$ is contained in a $3$-dimensional linear subspace of $\bP^5$, we are essentially done by theorem  \ref{thmsstable}. The detailed discussion is done in \S\ref{nonisolated1}. The remaining cases, $C$ is contained in a hyperplane or $C$ is linearly non-degenerate, are discussed in 
\S\ref{nonisolated2} and \S\ref{nonisolated3} respectively.  These non-degenerate cases produce the new strata $\epsilon$ and $\phi$.\end{proof}

\begin{remark}\label{remothers} We note that there are several results in literature that partially overlap with the results of theorem \ref{thmsingularities} (see especially Aluffi \cite{aluffising} and O'Grady \cite{ogrady}). 
\end{remark}

The general equations for the cubic fourfolds parametrized by the strata $\epsilon$ and $\phi$ of the theorem are:
\begin{eqnarray}\label{eqepsilon}
\epsilon:\ \ g(x_0,\dots,x_5)&=&\left| \begin{array}{ccc}
x_0&x_1&x_2+2a x_5\\
x_1&x_2-a x_5&x_3\\
x_2+2ax_5&x_3&x_4
\end{array}\right|\\
\notag&&+x_5^2l(x_0,\dots,x_4)+b x_5^3,
\end{eqnarray}
where $a,b\in \bC$ and $l$ is a linear form in $x_0,\dots,x_4$, and 
\begin{equation}\label{eqphi}
\phi:\ \ g(z_0,\dots,z_5)=\left| \begin{array}{ccc}
l_0&l_1&l_2\\
l_3&l_4&l_5\\
l_6&l_7&l_8
\end{array}\right|,
\end{equation} 
where $l_0,\dots,l_8$ are linear forms in $x_0,\dots,x_5$. A fourfold of type $\epsilon$ is a singular along a rational normal curve of degree $4$, which lies in the hyperplane $(x_5=0)$. An important invariant of a fourfold of type $\epsilon$ is the cross-ratio of the $4$ points cut on this rational curve by the linear form $l(x_0,\dots,x_4)$. The degenerate case (when the $4$ points are not distinct)  produces the surface $\sigma$:
\begin{equation}\label{eqsigma}
\sigma:\ \ g(x_0,\dots,x_5)=\left| \begin{array}{ccc}
x_0&x_1&x_2+2a x_5\\
x_1&x_2-a x_5&x_3\\
x_2+2ax_5&x_3&x_4
\end{array}\right|+cx_5^2 x_2+b x_5^3,
\end{equation}
for $a,b,c\in \bC$. Inside the surface $\sigma$, there are two curves that parameterize degenerate geometric situations. Namely, the curve $\tau$ (see (\ref{eqtau})) that was discussed in section \ref{sectminorbits}, and the curve $\chi$:
\begin{equation}\label{eqchi}
\chi:\ \ g(x_0,\dots,x_5)=\left| \begin{array}{ccc}
x_0&x_1&x_2+2a x_5\\
x_1&x_2-a x_5&x_3\\
x_2+2ax_5&x_3&x_4
\end{array}\right|+b x_5^3,
\end{equation}
for $a,b\in \bC$. The curves $\tau$ and $\chi$ intersect in the point $\omega$ (see (\ref{eqomega}) and remark \ref{defomega}), which  is in some sense the most degenerate point of the GIT quotient. The resulting adjacency diagram is given in figure \ref{gitboundary2}. For further details on the strata $\epsilon$ and $\phi$, we refer the reader to \S\ref{nonisolated2} and \S\ref{nonisolated3} respectively.

\begin{figure}[htb]
\begin{center}
$$
\xymatrix@R=.25cm{
&&&&{\ \ \ \ \ \ \ \ \ \ G^0_{\omega}\cong \SL(3)}\ar@{<-}[dd] \ar@{<-}[dddl]&\\
&{\ \ \ \ \ \ \ \ \ \ G^0_{\zeta}\cong T^4}\ar@{<-}[dd] \ar@{<-}[dddl]\ar@{<-}[ddrr]\\
&&&&{\ \ \ \ \ \ \ \ \ \ G^0_{\chi}\cong \SL(2)}\ar@{<-}[dd] \\
&{G^0_{\delta}}&&{G^0_{\tau}}\ar@{<-}[dl]\ar@{<-}[dr]&\\
{G^0_{\alpha}}&&{G^0_{\gamma}}&&{G^0_{\beta}}
}
$$
\end{center}
\caption{The stabilizers of general points in the boundary strata}\label{listr}
\end{figure}
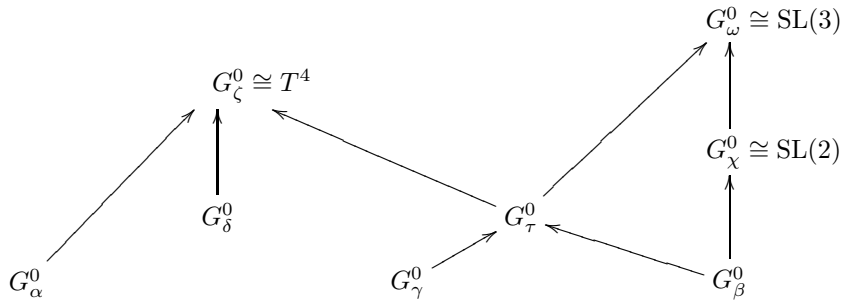

\begin{remark}
The boundary strata for a strict  GIT analysis are given in figure \ref{gitboundary}.  Figure \ref{gitboundary2} is obtained from figure \ref{gitboundary} by adding information about the possible singularities of a cubic fourfold. It is interesting to note that much of the extra structure of figure \ref{gitboundary2} can be also obtained by taking into account the possible stabilizers of cubic fourfolds  (compare figures \ref{gitboundary2} and \ref{listr}). The stratification given by stabilizers plays a key role in the analysis of the geometry of $\overline{\calM}$ from a GIT point of view (see the work of Kirwan \cite{kirwan,kirwanhyp}).
\end{remark}

\subsection{Strictly semistable cubic fourfolds with non-isolated singularities}\label{nonisolated1}
According to theorem \ref{thmsstable} a cubic fourfold containing a line, a conic or a rational elliptic curve of degree $4$ is not stable. Here, we prove that these are all the possibilities for a cubic fourfold to be singular along a curve contained in a $3$-dimensional linear subspace of $\bP^5$. 

\begin{proposition}\label{propless3}
Let $Y$ be a cubic fourfold with $\dim \Sing(Y)=1$. Assume that an irreducible $1$-dimensional component $C$ of the singular locus  is contained in a $3$-dimensional subspace of $\bP^5$. Then, $C$ is either  a line, a conic, or a component of the complete intersection of two quadrics in $\bP^3$. In particular, $Y$ is non-stable of type S1, S5, or S2 respectively. A general cubic $Y$ singular along a line, a conic, or an elliptic normal curve of degree $4$ is strictly semistable with the closure of its orbit containing a closed orbit of type $\alpha$, $\gamma$, or $\alpha$ respectively.
\end{proposition}
\begin{proof}
 Let $C\subset \Sing(Y)$ be a curve as in proposition, and denote by $\Pi:=\mathrm{Sec}(C)$ its secant variety. It is well known that a cubic contains the secant variety of its singular locus. Thus, $\Pi\subset Y$. Since the secant variety of a curve has always the expected dimension, it also follows that $\Pi$ is a linear subspace of $\bP^5$, coinciding with the linear span of $C$. Now the claim follows easily. For example, say $\dim \Pi=3$. Under appropriate choice of coordinates, we can assume $\Pi$ is given by $(x_4=x_5=0)$. It follows that the defining equation of $Y$ can be taken to be: 
$g(x_0,\dots,x_5)=x_4Q_1(x_0,\dots,x_5)+x_5Q_2(x_0,\dots,x_5)$,
for some quadrics $Q_1$ and $Q_2$. Let $q_1(x_0,\dots,x_3)$ and $q_2(x_0,\dots,x_3)$ be the quadrics obtained by restricting $Q_1$ and $Q_2$ to $\Pi$. If both $q_1$ and $q_2$ vanish, we get $\dim \Sing(Y)\ge2$. Otherwise, it is easily seen that  the intersection of  $\Pi$ with the singular locus $\Sing(Y)$ is the complete intersection of the quadrics $q_1$ and $q_2$. The conclusion follows.
\end{proof}

The stability  of cubic fourfolds with singular locus at least $2$-dimensional is settled by the following proposition:

\begin{proposition}\label{propsingsurf}
Let $Y$ be a cubic fourfold with $\dim \Sing(Y)\ge 2$. Then one of the following holds:
\begin{itemize}
\item[i)] the singular locus contains a plane or a quadric surface;
\item[ii)] the singular locus is the cone over the rational normal quartic, in which case $Y$ is the cone over the chordal cubic threefold (unique up to projective equivalence);
\item[iii)] the singular locus is the Veronese surface in $\bP^5$, in which case $Y$ is the secant variety of the Veronese surface (unique up to projective equivalence).
\end{itemize}
Furthermore, $Y$ is semistable if and only if it is the secant variety of the Veronese surface, in which case it gives the minimal orbit  $\omega$. 
\end{proposition}
\begin{proof}
The case of $\dim \Sing(Y)=3$ is easily settled ($Y$ is singular along a plane and unstable). Assume $\dim \Sing(Y)=3$ and let  $S$  be an irreducible surface in the singular locus. By cutting $Y$ with a generic hyperplane we obtain a cubic threefold singular along the irreducible curve $C$ cut on $S$ by this hyperplane. It is known that the degree of $C$ (and thus $S$)  is either $1$, $2$ or $4$ (cf. \cite[Prop. 4.2]{yokoyama}). For degree $1$ or $2$ the proposition follows immediately.

If $S$ is of degree $4$, we can assume that $S$ is non-degerate. Otherwise, one obtains a contradiction. It follows that $S$ is a surface of minimal degree in $\bP^5$, i.e. $S$ is:
\begin{itemize}
\item[-] a smooth rational quartic scroll,
\item[-] the Veronese surface,
\item[-] or the cone over the rational normal quartic curve in $\bP^4$. 
\end{itemize}
(see \cite[pg. 523]{griffithsharris}). Since $Y$ is singular along $S$, we have $\mathrm{Sec}(S)\subset Y$. The expected dimension of the secant variety of a surface is $5$.  For smooth surfaces, the expected dimension is the actual dimension with a single exception, the Veronese surface (see \cite[pg. 179]{griffithsharris}). Thus, if  $S$ is smooth, we must have that $S$ is the Veronese surface and $Y$ is the secant variety of $S$.  As discussed in \S\ref{sectminorbits}, the fourfold $Y$ gives the minimal orbit $\omega$. Finally, if $S$ is singular (i.e. the cone over the rational normal curve),  it is easily seen that $Y$ is the cone over the chordal cubic threefold.  It follows that $Y$ is unstable of type U6 (cf. lemma \ref{caseu6}). \end{proof}

\subsection{Cubic fourfolds singular along a rational normal curve of degree $4$}\label{nonisolated2} 
We now consider the case of cubic fourfolds singular along a curve that is contained (and is non-degenerate) in a hyperplane in $\bP^5$. We establish that a general fourfold of this type is stable, and that the corresponding locus in the moduli space of cubic fourfolds $\overline{\calM}$ is irreducible and $3$-dimensional. The key to this result is that we can find  a normal form for such fourfolds for which the projective invariants are clear.

\begin{proposition}
Let $Y$ be a cubic $4$-fold with $\dim \Sing(Y)=1$. Assume $C\subseteq \Sing(Y)$ is an irreducible curve whose linear span is $4$-dimensional. Then $C$ is a rational normal quartic in $\bP^4$. The equation of $Y$ can be taken of  the form (\ref{eqepsilon}). 
\end{proposition}
\begin{proof}
Let  $H$ be the hyperplane spanned by $C$ and  $X$  the cubic threefold obtained by restricting $Y$ to $H$. By construction, $X$ is singular along the curve $C$. Since $C$ is non-degenerate in $H$, it follows that $X$ is irreducible. By the analysis of cubic threefolds with non-isolated singularities (e.g. \cite[Prop. 4.2]{yokoyama}), it follows that $C$ is a rational normal curve.

To establish the  normal form (\ref{eqepsilon}) for $Y$, we note that  $X$ is the secant variety of $C$, i.e.    $X$ is the chordal cubic $3$-fold (see \cite[\S1]{allcock1}). Thus, we can choose coordinates such that the equation of $Y$ is given by 
\begin{equation}\label{temp1}
g(x_0,\dots,x_5)=F_0(x_0,\dots,x_4)+x_5Q(x_0,\dots,x_4)+x_5^2l(x_0,\dots,x_4)+b x_5^3,
\end{equation}
where $F_0$ is the equation of the chordal cubic and the hyperplane $H$ is given by $(x_5=0)$. The condition that $C$ is in the singular locus of $Y$ is equivalent to asking that the quadric given by $Q(x_0,\dots,x_4)$ in $H\cong \bP^4$ contains $C$. There exists a $6$-dimensional linear system of quadrics containing the rational normal $C$ (i.e. $\dim H^0(\bP^4,\calI_C(2))=6$). By using linear changes of coordinates, one sees that $Q$ is defined up to addition of elements from the Jacobian ideal of $F_0$. The partials of $F_0$ spans a $5$-dimensional subspace of $H^0(\bP^4,\calI_C(2))$. Thus, there is no loss of generality in assuming that in (\ref{temp1}) we have $Q=aQ_0$ for any quadric $Q_0$ containing $C$ and not lying in the Jacobian ideal. There is in fact a canonical choice for $Q_0$ (up to scalling). Namely, we choose $Q_0$ to be the unique quadric left invariant by $\SL(2)$ acting on $\bP^4\cong \Sym^4\bP^1$ via the action induced by the natural action on $\bP^1$ (see \cite[Ex. 10.12]{harris}). Concretely,  $Q_0(x_0,\dots,x_5)=4x_1x_3-3x_2^2-x_0x_4$ and  the normal form (\ref{eqepsilon}) is obtained immediately. 
\end{proof}

We have proved  that any cubic fourfold singular along a curve spanning a hyperplane in $\bP^5$ can be put in the normal form (\ref{eqepsilon}). There are two types of transformations preserving the equation (\ref{eqepsilon}): the scaling of the variable $x_5$ and the action of $\SL(2)$ on the  variables $x_0,\dots,x_4$ by means of the isomorphism $\bP^4\cong \Sym^4 \bP^1$ (induced by the rational normal curve of degree $4$). It follows that (up to projective transformations) the cubics of type $\epsilon$ depend on $3$ parameters. One of these parameters is quite geometric: it is the $j$-invariant of the $4$ points cut on the rational normal curve (the singular locus of the cubic) by the linear form $l$. Expanding on this observation, we obtain:

\begin{proposition}\label{caseepsilon} Let $Y$ be a cubic fourfold singular along a rational normal curve $C$ of degree $4$. Assume that defining equation of $Y$ is of the form (\ref{eqepsilon}). Then, one of the following holds:
\begin{itemize}
\item[i)] The hyperplane given by $l$ cuts $C$ in four distinct points. In this situation $Y$ is stable. The singularities of $Y$ along $C$ are of type $A_{\infty}$, except at the $4$ points where $l$ cuts $C$. At these points the singularities are of type $D_\infty$. 
\item[ii)] The hyperplane given by $l$ is simply tangent to $C$ in some point. Then, the orbit of $Y$ contains in its closure the orbit of a fourfold of type $\sigma$. A general point of $\sigma$ corresponds to a semi-stable cubic fourfold with closed orbit (otherwise, it degenerates to a fourfold of type $\tau$).
\item[iii)] The hyperplane given by $l$ intersects $C$ with multiplicity at least $3$, or $l$ is identically $0$. Then, the orbit of $Y$ contains in its closure the orbit of a fourfold of type $\chi$. For $a$ and $b$ not simultaneously $0$ in (\ref{eqchi}), a fourfold of type $\chi$ is semi-stable with closed orbit.  
\item[iv)] If $a=b=0$ and $l\equiv 0$ then $Y$ is the cone over the chordal cubic threefold and it is unstable. 
\end{itemize}
\end{proposition}
\begin{proof}
The stability statement follows by analyzing the possible singularities of $Y$, and concluding that, in the case i), $Y$ does not satisfy any of the conditions of theorem \ref{thmsstable}. If  the hyperplane given by $l$ becomes tangent to $C$, the fourfold $Y$ is no longer stable; it is of type S3 (with respect to the $1$-PS of weights $(2,1,0,-1,-2,0)$). It then follows that $Y$ degenerates to a fourfold of type $\sigma$. The case $\sigma$ is a particular case of $\beta$ and as such one can determine precisely when the orbit is closed (see the first part of lemma \ref{casebeta}). In particular, it follows that the surface $\sigma$ meets the other boundary components along the curve $\tau$. Similar arguments apply to the case $\chi$. 
\end{proof}

\begin{remark}\label{remchi}
Strictly speaking, one does not need to separate the case $\chi$ from the case $\sigma$. However, we choose to consider it as a separate case for two reasons. First, the stabilizer in the case $\chi$ is $\SL(2)$ versus $\bC^*$ for the generic point of $\sigma$. Secondly, the singularities along the rational normal curve are of transversal type $A_2$ for $\chi$ versus transversal type $A_1$ for a general point on $\sigma$.
\end{remark}

\subsection{Cubic fourfolds singular along a non-degenerate curve}\label{nonisolated3} 
In order to establish theorem \ref{thmsingularities}, the last case we have to consider is that of cubic fourfolds singular along a non-degenerate curve. We prove that in this case the cubic is determinantal (see equation (\ref{eqphi})) and generically stable. The first step for this is a geometric argument showing that there are only two possible types of non-degenerate curves that occur in the singular locus of a cubic fourfold.
\begin{proposition}\label{propprojection}
Let $Y$ be a cubic fourfold with $\dim \Sing(Y)=1$. Assume that $C\subset \Sing(Y)$ is an irreducible non-degenerate curve. Then $C$ is either a rational normal quintic or an elliptic normal sextic or degeneration of it.   
\end{proposition}
\begin{proof}
Let $C\subseteq \Sing(Y)$ be an irreducible non-degenerate curve of degree $d$. We choose a generic point $p\in C$ and project $Y$ onto a generic hyperplane $H\cong \bP^4$. Inside $H$ we obtain a sextic surface $S_1$ (a $(2,3)$ complete intersection) that parametrizes the lines passing through $p$.  The projection $C_1$ of $C$ is irreducible, non-degenerate, birational to $C$ (\cite[pg. 109]{arbarello}), and of degree $d-1$ (\cite[pg. 235]{harris}). Since $Y$ is singular along $C$, the surface $S_1$ contains $C_1$. It is easy to see that $S_1$ is singular along $C_1$. Note that the surface $S_1$ is non-degenerate and reduced (otherwise $\dim \Sing(Y)\ge 2$). For degree reasons, the only possibility for $S_1$ to be reducible is to be the union of two non-degenerate cubic scrolls. 

We repeat the construction with $C_1$ and $S_1$. By projecting from a generic point of $C_1$ we obtain a curve $C_2$ and a surface $S_2$. The curve $C_2$ is birational to $C$ and of degree $d-2$ in $\bP^3$. Since we project from a general point of $C_1$, and $S_1$ is a complete intersection, one checks that $S_1$ maps birationally onto its image. It is also clear that  $C_2\subset \Sing(S_2)$. Since we projected from a singular point of $S_1$, the surface $S_2$ will have degree $4$. We have two cases: either $S_2$ is irreducible or not. If  the quartic $S_2$ is irreducible, then the degree of $C_2$ is at most $3$ (a generic hyperplane section is an irreducible plane quartic, which has at most $3$ singular points). In this case, we conclude that $C_2$ is the twisted cubic. It follows that $C$ is a rational normal quintic in $\bP^5$. On the other hand, if $S_2$ is reducible, the only non-degenerate case is of $C_2$ being a component of a complete intersection of two quadrics. Thus $C_2$ has degree $3$ or $4$ and is rational or elliptic. The conclusion follows.
\end{proof}

We have the following dimension count for the space of cubic fourfolds singular along a given elliptic sextic curve or rational normal curve:
\begin{lemma}\label{computedim}
Let $C$ be an elliptic normal curve of degree $6$ (or rational normal curve of degree $5$) in $\bP^5$. Then, the linear system of cubics singular along $C$ is $1$-dimensional  (resp. $3$-dimensional). 
\end{lemma}
\begin{proof}
The cubics singular along $C$ correspond to the global sections of $\calI_C^2(3)$. Starting from the exact sequences:
$$0\to\calI_C\to\calO_{\bP^n}\to \calO_C\to 0$$
and 
$$0\to \calI_C^2\to \calI_C\to \calI_C/\calI_C^2\to 0,$$
one obtains (via a standard Riemann-Roch computation) the Euler characteristic: 
$$\chi(\calI_C^2(3))={n+3 \choose 3}-(2n-1)d+(n+2)(g-1),$$
where $n=5$ is the dimension of the projective space, $d$ and $g$ are the degree and genus of  $C$. In the rational case, we have the expected dimension $\chi(\calI_C^2(3))=4$. Similarly, if $C$ is an elliptic sextic, then $\chi(\calI_C^2(3))=2$.

To compute the actual dimensions of $H^0(\calI_C^2(3))$, we combine the Euler characteristic computation with some vanishing results. Specifically, by a result of Rathmann (see \cite[Thm. 1.1]{vermeire}) one is guaranteed the vanishing of  $H^i(\bP^n,\calI_C^2(k))=0$ for all $k\ge 3$ and $i>0$, provided that the degree of the curve is high with respect to the genus, e.g. $d\ge 2g+3$ suffices. In both cases considered here,  the vanishing result applies. Thus, $h^0(\calI_C^2(3))=\chi(\calI_C^2(3))$ and the lemma follows.
\end{proof}

The structure of the cubic fourfolds singular along an elliptic normal curve is described by the following result:
\begin{lemma}\label{detenc6}
Let $Y$ be  a cubic fourfold whose equation is given by the determinant of a $3\times 3$ matrix $A$ of linear forms (we call $Y$  determinantal). Then $Y$ is singular along the curve cut by the minors of  $A$. The singular locus of a general determinantal cubic fourfold is  precisely  an elliptic normal curve of degree $6$. Conversely, a cubic fourfold singular along an elliptic normal curve of degree $6$ is determinantal.
\end{lemma}
\begin{proof}
The fourfold $Y$ is determinantal if and only if it is a linear section $L$ of the secant variety (a cubic hypersurface in $\bP^8$) of the Segre fourfold (the image of $\bP^2\times \bP^2\hookrightarrow \bP^8$). In the generic case the singular locus of $Y$ is the intersection of the section $L$ with the Segre fourfold. The claim now follows from standard facts on the Segre embedding. Conversely, the fact that a cubic fourfold singular along an elliptic sextic is determinantal is a classical result, e.g. Room \cite[\S 7.11]{room} (see also \cite[Ex. 2.11]{grosspopescu}). \end{proof}

\begin{remark}\label{rempencil}
According to lemma \ref{computedim}, there exists a pencil of cubic fourfolds singular along a given elliptic normal curve $E\hookrightarrow \bP^5$. We note that this pencil is not trivial: there are $4$ special members of the pencil where the corresponding cubic is singular along a Veronese surface. Namely, the elliptic curve $E$ is embedded in $\bP^5$ by a complete linear system $|D|$ of degree $6$. Let $D'\in \mathrm{Pic}(E)$ such that $D=2D'$. It follows that the embedding of $E$ in $\bP^5$ can be factored as a composition of an embedding of $E$ in $\bP^2$ (given by $|D'|$) followed by a Veronese embedding $\bP^2\hookrightarrow \bP^5$. Thus, $E$ sits on a Veronese surface $S$. Then, $\Sec(S)$ is singular along $S$ (and thus $E$) and gives a section of $\calI_E^2(3)$. Since $D'$ is defined up to addition of points of order $2$ in the Jacobian of $E$, we obtain $4$ such sections.
\end{remark}

We now note that the rational curve case is  a specialization of the elliptic case.
\begin{lemma}\label{detrnc5}
A cubic fourfold singular along a rational normal curve of degree $5$ is determinantal. In particular, $Y$ is also singular along a line.
\end{lemma}
\begin{proof}
Let $C$ be a rational normal curve of degree $5$.  Without loss of generality, we can assume that $C$ is defined by the $2\times 2$ minors of the matrix:
$$M=\left(\begin{array}{cccc}
x_0&x_1&x_2&x_3\\
x_1&x_2&x_3&x_4\\
x_2&x_3&x_4&x_5
\end{array}\right)$$
Since the secant variety of $C$ is cut by the $3\times 3$ minors of $M$ (\cite[Prop. 9.7]{harris}), it follows that the four $3\times 3$ minors of $M$ give four linearly independent sections of $\calI_C^2(3)$. By lemma \ref{computedim},   we have $\dim H^0(\bP^5,\calI_C^2(3))=4$. Thus, any section of  $\calI_C^2(3)$ is a linear combination of the minors of $M$. In conclusion,  a cubic fourfold singular along a rational normal curve has the equation:
$$
g(x_0,\dots,x_5)=\left|\begin{array}{cccc}
a& b&c&d\\
x_0&x_1&x_2&x_3\\
x_1&x_2&x_3&x_4\\
x_2&x_3&x_4&x_5
\end{array}\right|
$$
for $a,b,c,d\in \bC$,  which can be then arranged  in a determinantal form. \end{proof}

We conclude:
\begin{proposition}\label{casephi}
Let $Y$ be an irreducible cubic fourfold singular locus along a non-degenerate curve $C$. Then $Y$ is determinantal (i.e. given by (\ref{eqphi})). A general determinantal cubic is stable. The locus of determinantal cubics in $\overline{\calM}$ forms an irreducible surface $\phi$, which meets the other components  of $\overline{\calM}\setminus \calM$ in the curve $\tau$.
\end{proposition}
\begin{proof}
We can restrict to the case $\dim \Sing(Y)=1$ (otherwise it is easily checked that $Y$ is determinantal). By proposition \ref{propprojection},  $Y$ is singular along a (non-degenerate) rational normal curve or an elliptic normal curve. According to lemmas   \ref{detenc6} and \ref{detrnc5}, in either of these cases $Y$ is determinantal. If $Y$ is general, then $Y$ is singular along a non-degenerate elliptic normal curve with  $A_{\infty}$ singularities. By theorem \ref{thmsstable}, it follows that $Y$ is stable. We define $\phi$ to be the closure in $\overline{\calM}$ of the locus of such fourfolds.  Clearly,  $\phi$ is $2$-dimensional: one dimension for the modulus of the elliptic curve, and one dimension for the choice of a point in the pencil of cubics singular along a given elliptic curve (see \ref{rempencil}).

The cubics parameterized by the curve $\tau$ (including the special points $\zeta$ and $\omega$) are determinantal.  It follows that $\tau$ is included in the surface $\phi$. The remaining part of the proposition follows from two observations. First, if $\Sing(Y)$ degenerates to a union of two or more curves,  then, for degree reasons, $Y$ is not stable (cf. theorem \ref{thmsstable}). Finally, if $Y$ is strictly semistable and determinantal, the minimal orbit contained in the orbit of $Y$ has to be of type $\tau$ (again for degree reasons).
\end{proof}

\section{Proof of the main results}\label{sectconclusion} 
At this point, the main theorems \ref{mainthm1} and \ref{mainthm2} are obtained as simple consequences of the results proved in the previous section.  For reader convenience, we collect the essential information about the boundary strata in table \ref{tableboundary}. We recall that the boundary strata $\alpha$--$\phi$ are closed. The various degeneracy loci inside these boundary strata (mostly corresponding to intersections of two or more strata) are given in table \ref{tableboundary2} (see also figure \ref{gitboundary2}). 

\begin{table}[htb]
\begin{center}
\renewcommand{\arraystretch}{1.25}
\begin{tabular}[2cm]{|c|c|c|c|c|}
\hline
Stratum & Equations&Discussion&Singularities & Dim. \\
\hline\hline
$\alpha$&(\ref{eqalpha})&\S\ref{sectminorbits}, esp. \ref{casealpha} &a  line and a quartic elliptic curve & 1\\
\hline
$\beta$&(\ref{eqbeta})&\S\ref{sectminorbits}, esp. \ref{casebeta} &two $\widetilde{E}_8$& 3\\
\hline
$\gamma$&(\ref{eqgamma})&\S\ref{sectminorbits}, esp. \ref{casegamma} & an $\widetilde{E}_7$ and a conic & 2\\
\hline
$\delta$&(\ref{eqdelta})&\S\ref{sectminorbits}, esp. \ref{casedelta} &three $\widetilde{E}_6$& 1\\
\hline
$\epsilon$&(\ref{eqepsilon})&\S\ref{nonisolated2}, esp. \ref{caseepsilon} & a rational quartic curve & 3\\
\hline
$\phi$&(\ref{eqphi})&\S\ref{nonisolated3}, esp. \ref{casephi} &an elliptic sextic curve & 2\\
\hline
\end{tabular}
\vspace{0.1cm}
\caption{The  boundary of $\calM\subset \overline{\calM}$}\label{tableboundary}
\end{center}
\end{table}

\begin{table}[htb]
\begin{center}
\renewcommand{\arraystretch}{1.25}
\begin{tabular}[2cm]{|c|c|c|c|c|}
\hline
Boundary & Equations&Specialization of & Dimension \\
\hline\hline
$\sigma$&(\ref{eqsigma})&$\beta$, $\epsilon$  &  2\\
\hline\hline
$\tau$&(\ref{eqtau})&  $\gamma$, $\phi$, $\sigma$ & 1\\
\hline
$\chi$&(\ref{eqchi})&$\sigma$ & 1\\
\hline\hline
$\omega$&(\ref{eqomega})&$\tau$, $\chi$ &  0\\
\hline
$\zeta$&(\ref{eqzeta})&$\alpha$, $\delta$, $\tau$&  0\\
\hline
\end{tabular}
\vspace{0.1cm}
\caption{The special strata of the boundary of $\calM\subset \overline{\calM}$}\label{tableboundary2}
\end{center}
\end{table}

\subsection{Proof of Theorem \ref{mainthm1}}
A cubic fourfold $Y$ is not properly stable if and only if it satisfies the degeneracy conditions of  theorem \ref{thmsstable}. The degeneracy conditions corresponding to cases S1, S2, and S5 are equivalent to saying that $Y$ is singular along a curve contained in a $3$-dimensional  linear subspace of $\bP^5$ (cf. proposition \ref{propless3}). The conditions S3, S4, and S6 are equivalent to requiring that $Y$ contains a singular point $p$ that deforms to a singularity of $\widetilde{E_r}$ for $r=8$, $7$, or $6$ respectively (see the discussion of section \ref{sectisolated}, especially theorem \ref{thmsimplesing}). We note that the proofs of section \ref{sectisolated} apply also to the case when the singularity at $p$ is non-isolated. \qed

\subsection{Proof of Theorem \ref{mainthm2}} According to theorem \ref{mainthm1} a cubic fourfold with simple singularities is stable.  In particular, it makes sense to talk about a moduli space $\calM$ of such fourfolds as a geometric quotient. Since the condition of simple singularities is open (a small deformation of A-D-E  singularities remains of type A-D-E) and $\SL(6)$-invariant, we obtain that $\calM$ is an open subset in $\calM^s$, which in turn is open in $\overline{\calM}$. We understand the boundary of $\calM$ in $\overline{\calM}$ as the union of $\overline{\calM}\setminus \calM^s$ and $\calM^s\setminus \calM$.

The boundary components of the GIT compactification $\calM^s\subset \overline{\calM}$ parametrize the minimal orbits of strictly semi-stable cubic fourfolds.  According to theorem \ref{thmsstable}, the strictly semistable cubic fourfolds with minimal orbits are of type $\alpha$--$\delta$, producing the boundary strata $\alpha$--$\delta$. The detailed analysis of these cases is done in  section \ref{sectminorbits}. Finally,  the points in $\calM^s\setminus \calM$ correspond to the orbits of stable cubic fourfolds with non-isolated singularities. By theorem \ref{thmsingularities}, the only relevant cases are those of the cubic fourfolds singular in a rational normal curve of degree $4$  or an elliptic normal curve of degree $6$. These two cases produce the boundary strata $\epsilon$ and $\phi$ (see \S\ref{nonisolated2} and \S\ref{nonisolated3}). For the discussion of the possible degeneracies and adjacencies we refer to the respective sections (see table \ref{tableboundary}).
\qed
\section{Further remarks}\label{sectcomments} 
\subsection{Relation to the stability of cubic threefolds} The analysis of stability for cubic fourfolds is a natural extension of the similar analysis for cubic threefolds of Allcock \cite{allcock1} and Yokoyama \cite{yokoyama}. Thus, there is no surprise that our results parallel to a large extent the results for cubic threefolds (compare for instance theorem \ref{thmsstable} to \cite[Thm. 1.3]{allcock1}). Here, we comment on a more direct link between the results for cubic threefolds and fourfolds.

 We recall that the construction of   Allcock--Carlson--Toledo \cite{allcock3fold}  associates to a cubic threefold $X\subset \bP^4$ the cubic fourfold $Y$ obtained as a triple cyclic cover of $\bP^4$ branched along $X$. It is not hard to see (due to the $\mu_3$ stabilizer) that the threefold $X$ gives a closed orbit if and only if the fourfold $Y$ does. This fact helped our analysis in two ways. First, the singularities allowed for stable cubic threefolds (i.e. $A_1,\dots,A_4$) are precisely  those that give stable singularities for the associated cubic fourfold $Y$. Thus, \cite[Theorem 1.1]{allcock1} anticipates our result on the stability of cubic fourfolds with simple singularities. Similarly, the boundary strata for cubic threefolds should correspond to boundary strata for cubic fourfolds. Indeed, this is the case: the strata $T$ and $\Delta$ of \cite[Thm. 1.2]{allcock1} correspond to our strata $\beta$ and $\delta$ respectively. Furthermore, the special point of $T$ corresponding to the chordal cubic corresponds to the special stratum $\chi\subset \beta$. Roughly speaking, the strata $\beta$ and $\delta$ are generated by $T$ and $\Delta$ in a natural way. For example,  $\delta$ can be interpreted as a $j$-line with $\Delta$  corresponding to the elliptic curve of $j$-invariant $0$. Thus, $\delta$ is obtained from $\Delta$ by allowing the $j$-invariant to vary.

\subsection{Relation to the stability of plane sextics} Hassett \cite[\S4.4]{hassett0} has noticed that the point $\omega$ corresponding to the secant variety of the Veronese surface plays a special role in the Hodge theoretical analysis of cubic fourfolds. Specifically, as one degenerates to $\omega$, one  obtains as the essential part of the limit mixed Hodge structure  the  Hodge structure of a degree $2$ $K3$ surface. Geometrically, this is explained by the following observation. Let $F_0$ be the equation of the secant to the Veronese surface $S\cong \bP^2$ in $\bP^5$. Consider a general pencil of cubics $(F_0+tF)_{t\in \bP^1}$ degenerating to $F_0$. Then, the base locus of the pencil cuts on the Veronese surface a plane sextic curve $C$. The degree $2$ $K3$ surface mentioned above is the double cover of $\bP^2$ along the sextic $C$ (the relation of Hodge structures is then obtained by applying the Clemens--Schmid exact sequence). 

By interpreting in a more intrinsic way the observation mentioned above, we find a natural relationship between the stability of cubic fourfolds and the stability of plane sextics (analyzed by Shah \cite{shah}). Namely, locally near $\omega$ (in the \'etale topology)  $\overline{\calM}$ is the quotient of an equivariant slice to the orbit $\omega$ by the stabilizer subgroup (Luna's slice Theorem \cite[Appendix D]{GIT}). The stabilizer corresponding to $\omega$ is a subgroup $H\cong \SL(3)$ of $G\cong \SL(6)$. The normal slice to the orbit of $\omega$ is computed as follows. Let  $W=H^0(\bP^5,\calO_{\bP^5}(1))\cong H^0(\bP^2,\calO_{\bP^2}(2))$ (isomorphism induced by the Veronese embedding). It follows that $W\cong \Sym^2 V$ as an $\SL(3)$-representation (where $V$ is the standard $\SL(3)$-representation) and then  
\begin{equation}\label{veroneserep}
\Sym^3(W)\cong \Sym^3(\Sym^2 V)\cong \bC\oplus \Gamma_{2,2}\oplus \Sym^6 V
\end{equation}
(cf. \cite[(13.15)]{fultonharris}). It is easy to see that the trivial summand in (\ref{veroneserep}) corresponds to the cubic fourfolds containing a fixed Veronese surface. Similarly, the $\Gamma_{2,2}$ summand  corresponds to the tangent directions to the orbit $\omega$. Thus, the $\SL(3)$-representation on the normal space is $\Sym^6 V$. We conclude that locally near $\omega$ the space $\overline{\calM}$ is isomorphic to $\Sym^6V/\SL(3)$, i.e. the affine cone over the GIT quotient for plane sextics. In particular, it follows that near $\omega$ the boundary structure of  $\overline{\calM}$ is essentially the same as that for plane sextics. Thus, we have the following natural correspondence between the boundary components adjacent to $\omega$ and the list of Shah for plane sextics (see \cite[Thm. 2.4]{shah}):
\begin{itemize}
\item[](Type \ II) $\beta$, $\gamma$, $\epsilon$, and $\phi$ correspond to II(1), II(2), II(3), and II(4) respectively;
\item[](Type III) $\sigma$, $\zeta$ correspond to III(1), and III(2) respectively;
\item[](Type IV) $\chi$ corresponds to IV;
\end{itemize}
(N.B. a geometric interpretation of this matching is easily seen by considering pencils of cubics degenerating to the secant to the Veronese, as  explained above).

\subsection{Stratification of the boundary}\label{secstratification} 
As mentioned above, the link between cubic fourfolds and degree $6$ plane curves is both at the level of GIT and Hodge theory. This motivates us to introduce a natural stratification of the boundary of the moduli space of cubic fourfolds. The stratification is defined in terms of the complexity (measured in ``types'') of the singularities for semistable cubic fourfolds. The stratification is meaningful from a GIT point of view, but, as in the case of Shah \cite{shah}, the real reason for introducing it is Hodge theoretic.

\begin{definition}\label{deftypeII}
We say $Y$ is a {\bf type I} cubic fourfold if $Y$ is smooth or has at worst simple isolated singularities. We say $Y$ is a {\bf type II} cubic fourfold if $Y$ is not of type I, and all the singularities of $Y$ are of the following types:
\begin{itemize}
\item[(0)] isolated simple singularities;
\item[(1)] isolated singularities of type $\widetilde{E_r}$ ($r=6,7,8$);
\item[(2)] non-isolated singularities of type $A_{\infty}$ (locally a double line);
\item[(3)] non-isolated singularities of type $D_{\infty}$ (simple pinch point).
\end{itemize} 
\end{definition}
Our results say that a cubic fourfold of type I is stable, and that one of type II is semi-stable. A generic boundary point belonging to a strata $\alpha,\dots,\phi$ is of type II. The boundary locus that does not parameterize cubic fourfolds of type II is precisely the surface $\sigma$. Thus, we define:
\begin{definition}\label{deftypeII}
Let $Y$ be a semi-stable cubic fourfold with minimal orbit. We say that $Y$ is of {\bf type III} if it is of type $\sigma$, but not of type $\chi$. We say that $Y$ is of {\bf type IV} if it is of type $\chi$.
\end{definition}

\begin{remark}
The type IV fourfolds are the ``most degenerate'' fourfolds that occur in the GIT analysis. This statement has meaning  both in terms of singularities and in terms of GIT (see \ref{defomega} and \ref{remchi}). A type IV has non-isolated singularities of transversal type $A_2$ (versus $A_1$ for type III). Also, the type IV fourfolds are the only semi-stable cubic fourfolds that are stabilized by a semisimple algebraic group. 
\end{remark}

In conclusion, we have stratified the boundary $\overline{\calM}\setminus \calM=\alpha\cup\dots \phi$ in 
\begin{itemize}
\item[] {\it type II boundary}: the open part $(\alpha\cup \phi)\setminus \sigma$,
\item[] {\it type III boundary}: the surface $\sigma\setminus \chi$,
\item[] and {\it type IV boundary}: the curve $\chi$.
\end{itemize}
This stratification is completely analogous to the stratification of Shah \cite{shah} for plane sextics.  The types correspond to  certain monodromy properties of the degenerations to a boundary point. Here, we only mention  a geometric peculiarity that plays an important role for the monodromy analysis: to each type II fourfold there is associated in a natural way an elliptic curve, or more precisely a $j$-invariant. Namely, a type II cubic is singular along an elliptic curve, along a rational curve with $4$ special points, or it has an $\widetilde{E}_r$ singularity (see the various lemmas where the boundary strata are analyzed); in each case the meaning of the $j$-invariant is clear.  The type III fourfolds  correspond to the case of $j$-invariant equal to $\infty$. For type IV fourfolds, one can not associate a meaningful $j$-invariant. Roughly speaking, the associated  elliptic curve becomes singular for the type III case, and completely vanishes for type IV.

\bibliography{references}
\end{document}